\documentclass[preprint,12pt]{elsarticle}

\usepackage{amssymb}
\usepackage{amsmath,amsfonts,mathrsfs,amssymb}
\usepackage{amsthm}

\numberwithin{equation}{section}

\begin{document}

\begin{frontmatter}



\title{Four dimensional hypersurfaces with proper mean curvature vector field in pseudo-Riemannian space forms}


\author[address]{Chao Yang}
\ead{yangch@nwnu.edu.com}

\author[address]{Jiancheng Liu}
\ead{liujc@nwnu.edu.cn}



\author[address2]{Li Du}
\ead{duli820210@cqut.edu.cn}


\address[address]{College of Mathematics and Statistics, Northwest Normal University, Lanzhou, China}
\address[address2]{School of Science, Chongqing University of Technology, Chongqing, China}

\begin{abstract}
  In this paper, we study four dimensional hypersurface $M^4_r$ with proper mean curvature vector field (i.e. $\Delta\vec{H}$ is proportional to $\vec{H}$) in pseudo-Riemannian space form $N^5_s(c)$, and show that it has constant mean curvature, and give the range of this constant. As an application, we get that biharmonic hypersurfaces in $N^5_s(c)$ are minimal in some specific cases, which partially confirms B.-Y. Chen's conjecture.
\end{abstract}


\end{frontmatter}


\section{Introduction}

Let $N^{n+1}_s(c)$ be a ($n+1$)-dimensional pseudo-Riemannian space form with index $0\leq s\leq n+1$ and constant sectional curvature $c$.
Specially, the pseudo-Riemannian space form $N^{n+1}_s(0)$ is isometric to ($n+1$)-dimensional pseudo-Euclidean space $\mathbb{E}^{n+1}_s$ with index $s$.
Let $x: M^n_r\to N^{n+1}_s(c)$ be an isometric immersion of
a pseudo-Riemannian hypersurface $M^n_r$ into $N^{n+1}_s(c)$.
Denote by $\vec{H}$ and $\Delta$ the mean curvature
vector field and the Laplace operator of $M^n_r$.
The hypersurface $M^n_r$ is said to have proper mean curvature vector field,
if it satisfies the equation
$$
\Delta\vec{H}=\lambda\vec{H},
$$
for some real constant $\lambda$.
Specially, when $\lambda=nc$, the hypersurface $M^n_r$ is biharmonic.

In 1988, B.-Y. Chen initiated the study of
hypersurface $M^n_r$ with proper mean curvature vector field in $\mathbb{E}^{n+1}_s$
in \cite{Chen 1988-1}, and proved that when $n=2, s=0$, surface $M^2$ is minimal, or an open part of a circular cylinder.
And then, A. Ferr\'{a}ndez and P. Lucus \cite{Ferrandez 1992-1} classified such non-minimal surfaces for $n=2, s=1$.
For $n=3$, the result of classification has not been gotten, but it was proved that hypersurface $M^3_r$ of $\mathbb{E}^4_s$ ($s=0, 1, 2$) has constant mean curvature ($s=0$ by F. Defever in \cite{Defever 1997} and T. Hasanis etc. in \cite{Hasanis 1995}, $s=1, 2$ by A. Arvanitoyeorgos etc. in \cite{Arvanitoyeorgos 2007, Arvanitoyeorgos 2009, Arvanitoyeorgos 2013}).

Naturally, there is a conjecture in \cite{Arvanitoyeorgos 2013} that: \emph{any hypersurface having proper mean curvature vector field in pseudo-Euclidean space $\mathbb{E}^{n+1}_s$ has constant mean curvature,}
which is closely related with the well known B.-Y. Chen's conjecture about biharmonic hypersurfaces in \cite{Chen 1991}: \emph{Any biharmonic hypersurface in $\mathbb{E}^{n+1}_s$ is minimal.} 
In \cite{Fu 2021}, Y. Fu, M.-C. Hong and X. Zhan illustrated the significance of solving Chen's conjecture for $n=4$, by analogy with the famous Bernstein problem.
In contrast to the Chen's conjecture, we realized that to solve the conjecture about hypersurfaces with proper mean curvature vector field for $n=4$ is necessary.

In 2021, Y. Fu and X. Zhan give an affirmative answer to this conjecture 
for $n=4, s=0$ in \cite{Fu 2021-2}. However, when $s>0$, the related research is difficult to carry out, since the shape operator of the hypersurface $M^4_r$ is not necessarily diagonalizable, and the principal curvatures may not all be real. In 2022, under the assumption that the hypersurface has constant scalar curvature and diagonalizable shape operator, L. Du and J. Ren obtained the same conclusion for $n=4, s>0$ in \cite{Du 2022}.
There are also some papers (\cite{Du 2017, Du 2023, Liu 2014, Liu 2016, Liu 2017}) studied the conjecture for $n>4$, but some strong restrictions are attached.

In this paper, we overcome the difficulties caused by the non-diagonalizable shape operator and the imaginary principal curvatures, and show that the hypersurface $M^4_r$ with proper mean curvature vector field in pseudo-Riemannian space form ${N}^{5}_s(c)$ has constant mean curvature in section 3 (cf. Theorem 3.1). 
Thus, the above conjecture about hypersurfaces with proper mean curvature vector field is true for $n=4$.

Once we know the mean curvature of the hypersurface $M^4_r$ is a constant, we continue in section 4 to estimate that constant.
When $M^4_r$ has at most two distinct principal curvatures, we get the value of the mean curvature $H$ (cf. Theorems 4.1, 4.6 and 4.8),
and find the value depends on $\varepsilon c$, $\varepsilon \lambda$ and the multiplies of the principal curvatures, where $\varepsilon=\langle\vec\xi,\vec\xi\rangle$, $\vec{\xi}$ denote a unit normal vector field to $M^4_r$.
Its a pity that when the number of distinct principal curvatures of $M^4_r$ is larger than 2, we do not obtain the value of the mean curvature. However, we give a value range of $H$ when the principal curvatures of $M^4_r$ are all real (cf. Theorem 4.12). As an application of these results, we can partially answer B.-Y. Chen's conjecture, details are as follows:
\begin{itemize}
  \item Biharmonic hypersurface $M^4_r$ in $N^5_s(c)$ with two distinct principal curvatures, which are imaginary, is minimal (cf. Corollary 4.11);
  \item Biharmonic hypersurface $M^4_r$ in $N^5_s(c)$ satisfying $c\varepsilon\leq 0$ and without imaginary principal curvatures, is minimal (cf. Corollary 4.13).
\end{itemize}

\section{Preliminaries}

\subsection{The formulas of hypersurface $M^n_r$ in $N^{n+1}_s(c)$}

Let $N^{n+1}_s(c)$ be a pseudo-Riemannian space form with index $s$ and constant sectional curvature $c$.
A non-zero vector $X$ in $N^{n+1}_s(c)$ is called \emph{time-like},
\emph{space-like} or \emph{light-like},
according to whether $\langle X, X\rangle$ is negative, positive or zero.

Let $M^n_r$ be a nondegenerate hypersurface in $N^{n+1}_s(c)$.
$\vec{\xi}$ denote a unit normal vector field to $M^n_r$,
then $\varepsilon=\langle\vec\xi,\vec\xi\rangle=\pm1$.
$\nabla$ and $\tilde{\nabla}$ denote the Levi-Civita connections of $M^n_r$
and $N^{n+1}_s(c)$, respectively. For any vector fields $X, Y$ tangent to $M^n_r$,
the Gauss formula is given by
\begin{equation*}
\widetilde\nabla_{X}Y=\nabla_{X}Y + h(X,Y)\vec{\xi},
\end{equation*}
where $h$ is the scalar-valued second fundamental form.
Denote by $A$, $\vec{H}$, and $H$ the shape operator of $M^n_r$ associated to
 $\vec{\xi}$, the mean curvature vector field and the mean curvature,
 then $\vec{H}=H\vec{\xi}$ and $H=\frac{1}{n}\varepsilon\text{tr}A$.

For any vector fields $X, Y, Z$ tangent to $M^n_r$, the Codazzi and Gauss equations are given by
\begin{equation}\label{Codazzi-td-eq}
\langle(\nabla_XA)Y, Z\rangle=
\langle(\nabla_YA)X, Z\rangle,
\end{equation}
and
\begin{equation*}
 R(X,Y)Z=c(\langle Y,Z\rangle X-\langle X,Z\rangle Y)+\varepsilon\langle A(Y),Z\rangle A(X)-\varepsilon\langle A(X),Z\rangle A(Y),
\end{equation*}
where $R(X,Y)Z=\nabla_{X}\nabla_{Y}Z-\nabla_{Y}\nabla_{X}Z-\nabla_{[X,Y]}Z$.

\subsection{The equivalent equations of $\Delta \vec{H}=\lambda \vec{H}$}

The equation $\Delta \vec{H}=\lambda \vec{H}$ can be rewritten as
$$
-\Delta\vec{H}-\text{trace}\,\tilde{R}(\text{d}x,\vec{H})\text{d}x=(\lambda-nc)\vec{H},
$$
which is equivalent to
$$
\tau_2(x)=(\lambda-nc)\tau(x),
$$
where $\tau(x)$ and $\tau_{2}(x)$ are the tension and bitension fields of $x$, respectively, and $\tilde{R}$ is the curvature tensor of $N^{n+1}_s(c)$.

Thus, according to \cite{Du 2017}, the hypersurface $M^n_r$ satisfies $\Delta \vec{H}=\lambda \vec{H}$ if and only if
\begin{equation}\label{proper-eq1}
 A(\nabla H)=-\frac{n}{2}\varepsilon H(\nabla H),
\end{equation}
and
\begin{equation}\label{proper-eq2}
\Delta H+\varepsilon H\text{tr}A^{2}=\lambda H,
\end{equation}
where the Laplace operator $\Delta$ acting on scalar-valued
function $f$ is given by
\begin{equation*}
\Delta f=-\sum^n_{i=1}\varepsilon_{i}(e_{i}e_{i}-\nabla_{e_{i}}e_{i})f,
\end{equation*}
with $\{e_i\}_{i=1}^n$ be a local orthonormal frame
such that $\langle e_i,e_i\rangle=\varepsilon_{i}=\pm1$.

\subsection{The shape operator of hypersurface $M^n_r$}

The tangent space $T_pM^n_r$ at $p\in M^n_r$ can be expressed as a direct sum of subspaces
$V_k$, $1\leq k\leq m$,
that are mutually orthogonal and invariant
under the shape operator $A$.
According to \cite[exercise 18, pp. 260-261]{Neill 1983}, there exists an integer $t$, with $0\leq t\leq m$, such that $A|_{V_i}$ (the restriction of $A$ on $V_i$), $1\leq i\leq t$,
has the form
\begin{equation*}
\begin{aligned}
A_i=\left(
  \begin{array}{ccccc}
     \lambda_{i}   &      &      &       &\\
        1   & \lambda_{i} &      &       &\\
        &       \ddots& \ddots   &       &\\
        &          &      1 &\lambda_{i} &\\
        &          &      &      1 & \lambda_{i}\\
  \end{array}
\right),
\end{aligned}
\end{equation*}
with respect to a basis $\mathfrak{B}_i=
\{u_{i_1}, u_{i_2}, \cdots, u_{i_{\alpha_i}}\}$ of $V_i$,
and $A|_{V_j}$, $t+1\leq j\leq m$, has the form
\begin{equation*}
\begin{aligned}
&\overline{A}_j=\left(
  \begin{array}{cccccccccc}
     \gamma_j&  \tau_j&          &         &           &         &       &         &         &\\
     -\tau_j& \gamma_j&         0&         &           &         &       &         &         &\\
      1&             0&  \gamma_j&  \tau_j &           &         &       &         &         &\\
      &              1&   -\tau_j& \gamma_j&          0&         &       &         &         &\\
      &               &         1&        0&   \gamma_j&   \tau_j&       &         &         &\\
      &               &          &        1&    -\tau_j& \gamma_j&      0&         &         &\\
      &               &          &         &     \ddots&   \ddots& \ddots&   \ddots&         &\\
      &               &          &         &           &        1&      0& \gamma_j&    \tau_j\\
      &               &          &         &           &         &      1&  -\tau_j&   \gamma_j\\
  \end{array}
\right), \tau_j\neq0,
\end{aligned}
\end{equation*}
with respect to a basis
$\overline{\mathfrak{B}}_j=\{u_{\bar{j}_1}, u_{\tilde{j}_1}, u_{\bar{j}_2}, u_{\tilde{j}_2}, \cdots, u_{\bar{j}_{\beta_j}},
u_{\tilde{j}_{\beta_j}}\}$ of $V_{j}$.
The inner products of the basis elements in $\mathfrak{B}_i$, $1\leq i\leq t$, and $\overline{\mathfrak{B}}_j$, $t+1\leq j\leq n$, are all zero except
\begin{equation*}
\langle u_{i_a}, u_{i_b}\rangle=\varepsilon_{i}=\pm1,\ \ a+b=\alpha_i+1,
\ \ 1\leq i\leq t,
\end{equation*}
and
\begin{equation*}
\langle u_{\bar{j}_c},u_{\bar{j}_d}\rangle=1=-\langle u_{\tilde{j}_c}, u_{\tilde{j}_d}\rangle,\ \ c+d=\beta_j+1,
\ \ t+1\leq j\leq m.
\end{equation*}
Certainly, the sum of the dimensions of $V_k$, $1\leq k\leq m$ is equal to $n$, i.e.
$$
\sum_{i=1}^t\alpha_i+2\sum_{j=t+1}^m\beta_j=n.
$$

Collecting all vectors in $\mathfrak{B}_1, \cdots, \mathfrak{B}_t,
\overline{\mathfrak{B}}_{t+1}, \cdots, \overline{\mathfrak{B}}_m$
in order, we get a basis
$$
\mathfrak{B}\!=\!\{u_{i_1}, u_{i_2}, \cdots, u_{i_{\alpha_i}}, u_{\bar{j}_1}, u_{\tilde{j}_1}, u_{\bar{j}_2}, u_{\tilde{j}_2}, \cdots, u_{\bar{j}_{\beta_j}},
u_{\tilde{j}_{\beta_j}}| 1\leq i\leq t, t+1\leq j\leq m\}
$$
of $T_xM^n_r$.
With respect to this basis $\mathfrak{B}$,
the shape operator $A$ of the hypersurface
$M^n_r$ in $N^{n+1}_s(c)$ can be expressed as an almost diagonal matrix
$$
A=\text{diag}\{A_{1}, \cdots, A_{t}, \overline{A}_{t+1}, \cdots, \overline{A}_{m}\}.
$$

Observe the form of $A_i$ and $\overline{A}_j$, with $1\leq i\leq t$, $t+1\leq j\leq m$,
we find $A_i$ has only a simple eigenvalue $\lambda_i$,
and $\overline{A}_j$ has eigenvalues $\gamma_j+\tau_j\sqrt{-1}$, $\gamma_j-\tau_j\sqrt{-1}$.
So, the shape operator $A$ has real eigenvalues
$$
\lambda_1, \lambda_2, \cdots, \lambda_t,
$$
and imaginary eigenvalues
$$
\gamma_{t+1}+\tau_{t+1}\sqrt{-1}, \gamma_{t+1}-\tau_{t+1}\sqrt{-1},
\cdots, \gamma_m+\tau_m\sqrt{-1}, \gamma_m-\tau_m\sqrt{-1}.
$$

According to the above explanations for the form of the shape operator, we conclude that for four dimensional hypersurface $M^4_r$, if there exists real principal curvatures, then the shape operator $A$ and corresponding metric matrix $G$ have possible forms:
\begin{scriptsize}
\begin{equation*}
\begin{array}{|c|c|c|c|}
 \hline \textbf{Form (I)} &  \textbf{Form (I\!I)} &  \textbf{Form (I\!I\!I)} \\
 \hline  \begin{array}{c}t=m=4,\\ \alpha_i=1, i=1, \cdots, 4.\end{array} &  \begin{array}{c}t=m=3,\\ \alpha_1=2, \alpha_2=\alpha_3=1.\end{array} & \begin{array}{c}t=m=2,\\ \alpha_1=\alpha_2=2.\end{array} \\
G=\left(
  \begin{array}{cccc}
     \varepsilon_{1}   &   0   &         0  &  0\\
        0   & \varepsilon_{2} &     0       & 0\\
      0  &   0       & \varepsilon_{3} & 0 \\
      0  &       0         &  0    & \varepsilon_{4}
  \end{array}
\right) &
G=\left(
  \begin{array}{cccc}
   0  &   \varepsilon_{1}     &     0       &        0\\
\varepsilon_{1}     & 0 &     0       &        0\\
     0    &      0             &   \varepsilon_{2}  &  0 \\
     0    &      0         &   0   & \varepsilon_{3}
  \end{array}
\right) &
G=\left(
  \begin{array}{cccc}
  0  &   \varepsilon_{1}     &      0       &   0\\
      \varepsilon_{1}     &   0  &     0     &   0\\
     0   &      0         &   0   &     \varepsilon_{2} \\
     0   &      0         &  \varepsilon_{2}    &  0
  \end{array}
\right)\\
  \mathfrak{B}=\{u_{1_1}, u_{2_1}, u_{3_1}, u_{4_1}\} &  \mathfrak{B}=\{u_{1_1}, u_{1_2}, u_{2_1}, u_{3_1}\}  &  \mathfrak{B}=\{u_{1_1}, u_{1_2}, u_{2_1}, u_{2_2}\}\\
  A=\left(\begin{array}{cccc}
     \lambda_{1}   &   0   &      0      &    0\\
       0     &    \lambda_{2}  &     0    &    0\\
      0    &           0      &  \lambda_{3}   &    0 \\
      0   &         0         &      0      &  \lambda_{4}
  \end{array}
\right) &  A=\left(\begin{array}{cccc}
     \lambda_{1}   &  0    &    0        & 0\\
         1  & \lambda_{1} &    0        & 0\\
      0   &       0        & \lambda_{2} &  0 \\
     0    &       0         &   0     &  \lambda_{3}
  \end{array}
\right) &
A=\left(\begin{array}{cccc}
     \lambda_{1}   &  0    &      0       &  0\\
          1 & \lambda_{1} &     0       & 0\\
    0    &        0        &   \lambda_{2} &  0 \\
    0     &          0       &     1   &  \lambda_{2}
  \end{array}
\right)
\\
 \hline  \textbf{Form (I\!V)} &  \textbf{Form (V)} &  \textbf{Form (V\!I)}\\
 \hline  \begin{array}{c}t=m=2,\\ \alpha_1=3, \alpha_2=1.\end{array} & \begin{array}{c}m=3, t=2,\\ \alpha_1=\alpha_2=1, \beta_3=1.\end{array} &  \begin{array}{c}m=2, t=1,\\ \alpha_1=2, \beta_2=1.\end{array}\\
   G=\left(
  \begin{array}{cccc}
   0  &  0  &  \varepsilon_{1}           &    0\\
   0   &  \varepsilon_{1}  &     0       &   0\\
     \varepsilon_{1}    &     0              &  0  & 0 \\
    0     &         0      &  0   &  \varepsilon_{2}
  \end{array}
\right) &
   G=\left(
  \begin{array}{cccc}
      \varepsilon_{1}   &   0   &         0   &  0 \\
        0   & \varepsilon_{2} &     0       & 0 \\
       0  &     0            &  1  &  0   \\
       0   &       0        &   0     &  -1
  \end{array}
\right) &
G=\left(
  \begin{array}{cccc}
    0 &  \varepsilon_{1}    &    0        &  0 \\
        \varepsilon_{1}      &  0 &      0       & 0\\
      0   &     0            &  1 &      0\\
      0   &     0         &   0   & -1
  \end{array}
\right)\\
\mathfrak{B}=\{u_{1_1}, u_{1_2}, u_{1_3}, u_{2_1}\} &  \mathfrak{B}=\{u_{1_1}, u_{2_1}, u_{\bar{3}_1}, u_{\tilde{3}_1}\} & \mathfrak{B}=\{u_{1_1}, u_{1_2}, u_{\bar{2}_1}, u_{\tilde{2}_1}\}\\
 A=\left(\begin{array}{cccc}
     \lambda_{1}   &   0   &    0        &  0\\
        1   & \lambda_{1} &      0      &  0\\
      0  &           1      &  \lambda_{1}  & 0 \\
      0  &         0       &    0  & \lambda_{2}
  \end{array}
\right) &
 A=\left(\begin{array}{cccc}
      \lambda_{1}   &    0    &       0       &  0  \\
         0     &   \lambda_{2}  &       0      &  0  \\
        0   &             0      & \gamma_{3} &   \tau_3   \\
        0   &         0         &    -\tau_3    &  \gamma_{3}
  \end{array}
\right) &  A=\left(\begin{array}{cccc}
    \lambda_{1}  &   0     &     0         &   0  \\
           1   &  \lambda_{1}  &     0         &  0\\
        0   &         0         &   \gamma_{2}  &   \tau_2\\
       0    &        0          &   -\tau_2    &  \gamma_{2}
  \end{array}
\right)
\\
\hline
\end{array}.
\end{equation*}
\end{scriptsize}
And if the principal curvatures of $M^4_r$ are all imaginary, then $A$ and $G$ have possible forms:
\begin{scriptsize}
\begin{equation*}
\begin{array}{|c|c|}
\hline \textbf{Form (V\!I\!I)} &  \textbf{Form (V\!I\!I\!I)}\! \\
\hline  \begin{array}{c}m=1, t=0,\\ \beta_1=2.\end{array} & \begin{array}{c}m=2, t=0\!\\ \beta_1=\beta_2=1.\end{array}\!\\
G=\!\left[\!
  \begin{array}{cccc}
  \! 0 \! &  \!0  \!  &  \!    1  \!    &\!0 \!\\
      \!0 \!    &\!0 \!&   \!  0   \!    &\!-1 \!\\
     \! 1  \! &     \! 0          \! &\!0\!& \!0  \!\\
    \!  0 \!  &  \!     -1    \!     & \!  0 \!  & \!0 \!\\
 \end{array}
\!\right] \!&\!
G=\!\left[\!
  \begin{array}{cccc}
    \!1\! &\! 0\!&\! 0\!           &\!0 \!\\
   \! 0 \! & \!-1\! &   \! 0      \!  &\!0 \!\\
    \!  0 \!  &  \!  0   \!          &\!1\!&\! 0 \! \\
     \!0  \!  &   \!       0  \!    & \!0 \! &\! -1 \!\\
 \end{array}\!
\right]\!\\
\mathfrak{B}\!=\!\{u_{\bar{1}_1}, u_{\tilde{1}_1}, u_{\bar{1}_2}, u_{\tilde{1}_2}\}\!& \!\mathfrak{B}\!=\!\{u_{\bar{1}}, u_{\tilde{1}}, u_{\bar{2}}, u_{\tilde{2}}\}\!\\
A\!=\!\left[\!\begin{array}{cccc}
     \! \gamma_{1}\!& \! \tau_1 \!   &   \!   0  \!     &  \!0 \!\\
      \!    -\tau_1\!   &\! \gamma_{1}\! &  \!   0   \!    &\! 0 \!\\
    \! 1  \!  &  \!      0   \!      &\!\gamma_{1}\!& \! \tau_1  \!\\
   \! 0  \!  &      \!   1  \!    &  \!  -\tau_1\!   &\! \gamma_{1} \!\\
  \end{array}
\!\right]\! &\!
A\!=\!\left[\!\begin{array}{cccc}
    \! \gamma_{1}\!& \! \tau_1 \!  &  \!  0   \!     & \! 0 \!\\
     \!   -\tau_1\!   &\! \gamma_{1}\! &  \!    0   \!   & \! 0 \!\\
     \! 0 \! &  \!         0  \!    &\!\gamma_{2}\!& \! \tau_2  \! \\
    \!  0 \! &   \!      0    \!   &  \!  -\tau_2\!   &\! \gamma_{2} \!\\
  \end{array}\!
\right]\!
\\
\hline
\end{array}.
\end{equation*}
\end{scriptsize}
From these forms of $A$, we can calculate $\text{tr}A$ and $\text{tr}A^2$, combining $\text{tr}A=4\varepsilon H$, we find
\begin{equation}\label{trA--td---sec2.3}
\begin{cases}
\text{tr}A=\sum_{i=1}^t\alpha_i\lambda_i+2\sum_{j=t+1}^m\beta_j\gamma_j=4\varepsilon H,\\
\text{tr}A^2=\sum_{i=1}^t\alpha_i\lambda_i^2+2\sum_{j=t+1}^m\beta_j(\gamma_j^2-\tau_j^2).
\end{cases}
\end{equation}

\section{The result that $M^4_r$ has constant mean curvature}

\vskip.2cm
\noindent
{\bf Theorem 3.1}\quad \emph{Let $N^{5}_s(c)$ be a $5$-dimensional pseudo-Riemannian space form with constant sectional curvature $c$, and $M^4_r$ be a nondegenerate hypersurface of $N^{5}_s(c)$ with proper mean curvature vector field,
then $M^4_r$ has constant mean curvature.}



To prove Theorem 3.1, we adopt contradiction method. Assume that $H$ is not a constant, we need to deduce a contradiction. When the principal curvatures of $M^4_r$ are all imaginary, \eqref{proper-eq1} implies $-2\varepsilon H$ (real) is an eigenvalue of $A$, a contradiction.
When $M^4_r$ has real principal curvatures, there are six possible forms of $A$ (see section 2).
These different forms of $A$ will generate different expressions about connection coefficients.
 To derive contradictions, we need to use these expressions. So, we cannot deal with these six forms uniformly.

\subsection{The shape operator has the form (I)}

\vskip.2cm
\noindent
{\bf Proposition 3.2}\quad \emph{Let $M^4_r$ be a nondegenerate hypersurface of $N^{5}_s(c)$ with proper mean curvature vector field. Suppose that the shape operator $A$ of $M^4_r$ has the form (I),
then $M^4_r$ has constant mean curvature.}

\vskip.2cm

When $M^4_r$ has at most three distinct principal curvatures, the proof has been provided by L. Du etc. in \cite{Du 2017}.
So, we suppose the principal curvatures $\lambda_1, \lambda_2, \lambda_3$ and $\lambda_4$ of $M^4_r$ are distinct with each other.
At the beginning, we give some equations and Lemmas, under the assumption that the mean curvature $H$ of $M^n_r$ is not a constant, and $\nabla H$ is in the direction $u_1$.


\subsubsection{Some equations and Lemmas}

\vskip.2cm
Denote $u_i=u_{i_1}$, with $1\leq i\leq 4$. Let $\nabla_{u_i}u_j=\Gamma^k_{ij}u_k$, $i, j=1, 2, 3, 4$.
we get from compatibility of the connection $\nabla$ that
\begin{equation}\label{5---compatibility-eq1--4}
\Gamma_{ki}^{j}=-\varepsilon_i\varepsilon_j\Gamma_{kj}^{i},
\end{equation}
with $i, j, k=1, 2, 3, 4$.
Assume that $H$ is not a constant, and $\nabla H$ is in the direction of $u_1$, then
\eqref{proper-eq1} implies $\lambda_1=-2\varepsilon H$, and
we have
\begin{equation}\label{5---B(H)--4case2}
u_{1}(H)\neq0,\ u_{2}(H)=u_{3}(H)=u_{4}(H)=0.
\end{equation}
As the equation \eqref{5---B(H)--4case2}, we obtain from $(\nabla_{u_i}u_j-\nabla_{u_j}u_i)(H)=[u_i, u_j](H)$ that
\begin{equation}\label{5---symmetry--4case2}
\Gamma^{1}_{ij}=\Gamma^{1}_{ji},\quad i, j=2, 3, 4.
\end{equation}

Combining the equations \eqref{5---compatibility-eq1--4}, \eqref{5---B(H)--4case2} and \eqref{5---symmetry--4case2},
we deduce from Codazzi equation that
\begin{equation}\label{5---codazzi-eq--1}
\Gamma_{1i}^{1}=0,\ \ \Gamma_{ij}^{1}=0,\ \text{with}\ i\neq j,\ i, j=2, 3, 4,
\end{equation}
and
\begin{equation}\label{5---codazzi-eq--2}
\begin{cases}
u_i(\lambda_j)=(\lambda_i-\lambda_j)\Gamma_{ji}^{j},\ i\neq j,\\
(\lambda_i-\lambda_j)\Gamma_{ki}^{j}=(\lambda_k-\lambda_j)\Gamma_{ik}^{j},\ i, k\neq j.
\end{cases}
\end{equation}
Using Gauss equation for $\langle R(u_1, u_i)u_1, u_i\rangle$, $i=2, 3, 4$,
combining \eqref{5---compatibility-eq1--4}, \eqref{5---symmetry--4case2} and \eqref{5---codazzi-eq--1}, we have
\begin{equation}\label{5---Gauss-eq--1}
\begin{aligned}
u_1(\Gamma_{i1}^i)=-(\Gamma_{i1}^i)^2+(2H\lambda_i-c)\varepsilon_1,\ i=2, 3, 4.
\end{aligned}
\end{equation}

Considering the expressions \eqref{5---compatibility-eq1--4} and \eqref{5---B(H)--4case2}, the equation (\ref{proper-eq2}) can be written as
\begin{equation}\label{5--proper-td-1}
u_1u_1(H)+\sum_{i=2}^4\Gamma_{i1}^{i}u_1(H)-\varepsilon \varepsilon_1H\text{tr}A^2+\varepsilon_1\lambda H=0.
\end{equation}

In the following, we discuss $
f_k:=\sum_{i=2}^4(\Gamma_{i1}^i)^k
$, with $k=1, 2, \cdots, 5$, and get the following Lemmas 3.3 and 3.4. And then, applying these two Lemmas, we prove $u_i(\Gamma_{j1}^j)=u_i(\lambda_j)=0$, with $i, j=2, 3, 4$, i.e. Lemma 3.5,
which will play an important role in the proof of Proposition 3.2.

\vskip.2cm
\noindent
{\bf Lemma 3.3}\quad \emph{We have
\begin{equation}\label{5---lemma5.2--1}
\begin{aligned}
f_2=&-u_1(f_1)+12H^2\varepsilon\varepsilon_1-3c\varepsilon_1;\\
f_3=&\frac{1}{2}u_1^{(2)}(f_1)-(4\varepsilon H^2+c)\varepsilon_1f_1-24\varepsilon\varepsilon_1 Hu_1(H);\\
f_4=&-\frac{1}{6}u_1^{(3)}(f_1)+\frac{4}{3}(4\varepsilon H^2+c)\varepsilon_1u_1(f_1)+\frac{20}{3}\varepsilon\varepsilon_1Hu_1(H)f_1
+8\varepsilon\varepsilon_1(u_1(H))^2\\
&+16\varepsilon\varepsilon_1 Hu_1^{(2)}(H)-32H^4+2\varepsilon H^2\lambda-12\varepsilon H^2c+3c^2;\\
f_5=&\frac{1}{24}u_1^{(4)}(f_1)-\frac{5}{6}(4\varepsilon H^2+c)\varepsilon_1u_1^{(2)}(f_1)
-\frac{25}{3}\varepsilon \varepsilon_1 Hu_1(H)u_1(f_1)-[16H^4\\
&-\frac{13}{3}Hu_1^{(2)}(H)-\frac{1}{3}\varepsilon \varepsilon_1(u_1(H))^2+8\varepsilon H^2c+c^2]f_1-8\varepsilon\varepsilon_1 Hu_1^{(3)}(H)\\
&-\frac{20}{3}\varepsilon\varepsilon_1 u_1(H)u_1^{(2)}(H)+[16(8+\frac{2}{3}\varepsilon\varepsilon_1)H^3+38\varepsilon Hc
-\frac{5}{3}H\varepsilon\lambda]u_1(H),
\end{aligned}
\end{equation}
where $u_1^{(k)}(f_1)$ and $u_1^{(k)}(H)$ with $k=2, 3, 4$ denote the $k$-order derivatives of $f_1$ and $H$ along $u_1$, respectively. }

\vskip.2cm
{\bf Proof}\quad
Since $\sum_{i=1}^4\lambda_i=\text{tr}A=4\varepsilon H$ and $\lambda_1=-2\varepsilon H$, we know
$$
\sum_{i=2}^4\lambda_i=6\varepsilon H.
$$
Taking sum for $j$ from $2$ to $4$ in the first equation of \eqref{5---codazzi-eq--2} with $i=1$, we find
\begin{equation}\label{5---lemma5.2-proof-2}
\begin{aligned}
\sum_{i=2}^4\lambda_i\Gamma_{i1}^i=-2\varepsilon Hf_1-6\varepsilon u_1(H).
\end{aligned}
\end{equation}
Multiply $\lambda_j$ on both sides of the first equation in \eqref{5---codazzi-eq--2} with $i=1$, and then we know
\begin{equation*}
\begin{aligned}
\sum_{i=2}^4\lambda_i^2\Gamma_{i1}^i=-2\varepsilon H\sum_{i=2}^4\lambda_i\Gamma_{i1}^i-\frac{1}{2}u_1(\text{tr}A^2)+4Hu_1(H).
\end{aligned}
\end{equation*}
Differentiate $\sum_{i=2}^4\lambda_i(\Gamma_{i1}^i)^k$ along $u_1$, combining \eqref{5---codazzi-eq--2} and \eqref{5---Gauss-eq--1},
we get
\begin{equation*}
\begin{aligned}
(k+1)\sum_{i=2}^4\lambda_i(\Gamma_{i1}^i)^{k+1}=&-u_1(\sum_{i=2}^4\lambda_i(\Gamma_{i1}^i)^k)+2k\varepsilon_1H\sum_{i=2}^4\lambda_i^2(\Gamma_{i1}^i)^{k-1}\\
&-k c\varepsilon_1\sum_{i=2}^4\lambda_i(\Gamma_{i1}^i)^{k-1}-2\varepsilon Hf_{k+1},\ \ k=1, 2,
\end{aligned}
\end{equation*}
which together with \eqref{5---lemma5.2-proof-2} and the above equation gives the expressions of $\sum_{i=2}^4\lambda_i(\Gamma_{i1}^i)^k$, with $k=2, 3$.
Take the sum of $i$ from 2 to 4 in \eqref{5---Gauss-eq--1}, we can express $f_2$ as
\begin{equation*}
f_2=-u_1(f_1)+12H^2\varepsilon\varepsilon_1-3c\varepsilon_1.
\end{equation*}
Multiplying $(\Gamma_{i1}^i)^k$ on both sides of \eqref{5---Gauss-eq--1}, and then taking sum for $i$, we obtain
\begin{equation*}
\begin{aligned}
f_{k+2}=-\frac{1}{k+1}u_1(f_{k+1})+2\varepsilon_1 H\sum_{i=2}^4\lambda_i(\Gamma_{i1}^i)^{k}-c\varepsilon_1f_{k},\ k=1, 2, 3.
\end{aligned}
\end{equation*}
As the expressions of $\sum_{i=2}^4\lambda_i(\Gamma_{i1}^i)^k$, with $k=1, 2, 3$, we conclude from the above equation
that \eqref{5---lemma5.2--1} holds.
$\hfill\square$

\vskip.2cm
\noindent
{\bf Lemma 3.4}\quad \emph{For $i=2, 3, 4$, we have
$$
u_i(f_1)=0.
$$}

{\bf Proof}\quad With the notions $f_k$, $k=1, 2, \cdots, 5$, we find the equations
\begin{equation}\label{5---fkrelations}
\begin{cases}
f_1^4-6f_1^2f_2+3f_2^2+8f_1f_3-6f_4=0,\\
f_1^5-5f_1^3f_2+5f_1^2f_3+5f_2f_3-6f_5=0
\end{cases}
\end{equation}
(cf. \cite{Fu 2021}) also hold. Substitute the expressions of $f_k$, $k=1, 2, \cdots, 5$ in Lemma 3.3 into \eqref{5---fkrelations}, we have
\begin{equation*}
\begin{aligned}
F_1:=u_1^{(3)}(f_1)+4f_1u_1^{(2)}(f_1)+3(u_1(f_1))^2+[6f_1^2-104\varepsilon\varepsilon_1H^2+10\varepsilon_1c]u_1(f_1)\\
+f_1^4-104\varepsilon\varepsilon_1H^2f_1^2+10\varepsilon_1cf_1^2-232\varepsilon\varepsilon_1Hu_1(H)f_1-48\varepsilon\varepsilon_1(u_1(H))^2\\
-96\varepsilon\varepsilon_1Hu_1^{(2)}(H)+624H^4-144\varepsilon H^2c-12\varepsilon H^2\lambda+9c^2=0,
\end{aligned}
\end{equation*}
and
\begin{equation*}
\begin{aligned}
F_2:=&-u_1^{(4)}(f_1)+(-10u_1(f_1)+10f_1^2+200\varepsilon\varepsilon_1H^2-10c\varepsilon_1)u_1^{(2)}(f_1)+4f_1^5\\
&+[20f_1^3
+20(4\varepsilon H^2+c)\varepsilon_1f_1+680\varepsilon\varepsilon_1Hu_1(H)]u_1(f_1)-320\varepsilon\varepsilon_1H^2f_1^3\\&+40\varepsilon_1f_1^3c
-480\varepsilon\varepsilon_1Hu_1(H)f_1^2+104Hu_1^{(2)}(H)f_1
+192\varepsilon\varepsilon_1Hu_1^{(3)}(H)\\&-192\varepsilon H^2cf_1-1344H^4f_1+8\varepsilon\varepsilon_1(u_1(H))^2f_1+160\varepsilon\varepsilon_1u_1(H)u_1^{(2)}(H)\\&+36c^2f_1-128(69+2\varepsilon\varepsilon_1)H^3u_1(H)
+(528c+40\lambda)\varepsilon Hu_1(H)=0.
\end{aligned}
\end{equation*}

Subsequently, we will eliminate $u_1^{(k)}(f_1)$, $k=1, 2, 3, 4$ from the equations $F_1=0$ and $F_2=0$. Let $F_3=\frac{1}{4}(u_1(F_1)+F_2)-f_1F_1$, then $F_3=0$, i.e.
\begin{equation}\label{5---lemma5.3-proof--4}
\begin{aligned}
&12\varepsilon\varepsilon_1H^2u_1^{(2)}(f_1)+[36\varepsilon\varepsilon_1H^2f_1+30\varepsilon\varepsilon_1Hu_1(H)]u_1(f_1)
+12\varepsilon\varepsilon_1H^2f_1^3\\&+30\varepsilon\varepsilon_1Hu_1(H)f_1^2+(19\varepsilon\varepsilon_1+13)Hu_1^{(2)}(H)f_1
-4\varepsilon\varepsilon_1(u_1(H))^2f_1\\&+6(8c+\lambda)\varepsilon H^2f_1
+12\varepsilon\varepsilon_1Hu_1^{(3)}(H)-4\varepsilon\varepsilon_1u_1(H)u_1^{(2)}(H)\\&-480H^4f_1-(792
+32\varepsilon\varepsilon_1)H^3u_1(H)+(30c+2\lambda)\varepsilon Hu_1(H)=0.
\end{aligned}
\end{equation}
Take $F_4=2H(u_1(F_3)-12\varepsilon\varepsilon_1H^2F_1+f_1F_3)-9u_1(H)F_3$,
then $F_4=0$, i.e.
\begin{equation}\label{5---lemma5.3-proof--6}
\begin{aligned}
a_1u_1(f_1)+a_1f_1^2+a_2f_1
+a_3=0,
\end{aligned}
\end{equation}
where
\begin{equation*}
\begin{aligned}
a_1=&(294\varepsilon\varepsilon_1+78)Hu_1^{(2)}(H)-654\varepsilon\varepsilon_1(u_1(H))^2+4608H^4
+36(\lambda-12c) \varepsilon H^2,\\
a_2=&(186\varepsilon\varepsilon_1+78)H^2u_1^{(3)}(H)-(471\varepsilon\varepsilon_1+273)Hu_1(H)u_1^{(2)}(H)\\&
+108\varepsilon\varepsilon_1(u_1(H))^3+[(13392-192\varepsilon\varepsilon_1)H^2
-6(90c+13\lambda)\varepsilon] H^2u_1(H),\\
a_3=&72\varepsilon\varepsilon_1H^2u_1^{(4)}(H)-276\varepsilon\varepsilon_1Hu_1(H)u_1^{(3)}(H)+108\varepsilon\varepsilon_1(u_1(H))^2u_1^{(2)}(H)
\\
&-24\varepsilon\varepsilon_1H(u_1^{(2)}(H))^2+[(2160-192\varepsilon\varepsilon_1)H^2+(180c+12\lambda)\varepsilon] H^2u_1^{(2)}(H)\\&+[72(147+4\varepsilon\varepsilon_1)H^2-(630c+42\lambda)\varepsilon] H(u_1(H))^2-44928\varepsilon\varepsilon_1H^7\\&
+(10368c+864\lambda)\varepsilon_1H^5-648\varepsilon\varepsilon_1H^3c^2.
\end{aligned}
\end{equation*}

By acting on \eqref{5---lemma5.3-proof--6} with $u_1$, and then combining \eqref{5---lemma5.3-proof--4}, we have
\begin{equation*}
(b_1f_1+b_2)u_1(f_1)+b_1f_1^3+b_3f_1^2+b_4f_1+b_5=0,
\end{equation*}
where
{\small \begin{equation*}
\begin{aligned}
b_1=&-768 H^4\varepsilon\varepsilon_1 + 6(12 c - 6\lambda) H^2\varepsilon_1 -( 13 \varepsilon\varepsilon_1+ 49) H u_1^{(2)}(H) + 109 (u_1(H))^2,\\
b_2=&(8304\varepsilon\varepsilon_1- 64) H^4 u_1(H) + (52\varepsilon\varepsilon_1+ 160) u_1^{(3)}(H) H^2+ 363 (u_1(H))^3\\
 &- (104 \varepsilon\varepsilon_1+ 642) H u_1(H) u_1^{(2)}(H) - (252 c +20\lambda) H^2 u_1(H)\varepsilon_1,\\
b_3=&3840 H^4 u_1(H)\varepsilon\varepsilon_1+(26\varepsilon\varepsilon_1+98)u_1^{(3)}(H) H^2 + 6 \varepsilon_1 (\lambda- 12 c)H^2 u_1(H) \\& + 327  (u_1(H))^3
- (13 \varepsilon\varepsilon_1 + 485) H u_1(H) u_1^{(2)}(H),\\
b_4=&12(13\varepsilon\varepsilon_1+31) u_1^{(4)}(H) H^3- 2(715 \varepsilon\varepsilon_1+1021) H^2 (u_1^{(2)}(H))^2+ 368640 H^8\\&+ 96(372\varepsilon\varepsilon_1 - 43) H^5 u_1^{(2)}(H) +(923\varepsilon\varepsilon_1 + 1973) H (u_1(H))^2 u_1^{(2)}(H)\\& - 18(13 \varepsilon\varepsilon_1+11) u_1^{(3)}(H) H^2 u_1(H)+ 6(464c+ 61\lambda)\varepsilon_1 H^2 (u_1(H))^2\\& + 36(96 c^2 -  \lambda^2 + 4 c\lambda) H^4+ 96(603 \varepsilon\varepsilon_1- 16) H^4 (u_1(H))^2- 436 (u_1(H))^4\\& - 576\varepsilon (124 c + 3\lambda)H^6+ 12(26 c\varepsilon -172 c\varepsilon_1 -13 \lambda\varepsilon -47\lambda\varepsilon_1) H^3 u_1^{(2)}(H),\\
b_5=&-12(144 c^2+ \lambda^2 +3  c \lambda) H^3 u_1(H) - 12(13\varepsilon\varepsilon_1+103) u_1^{(3)}(H) H^2 u_1^{(2)}(H)\\&+ 144 u_1^{(5)}(H) H^3 - 2(195 c\varepsilon   + 13\lambda\varepsilon+ 1779 c\varepsilon_1 +97 \lambda\varepsilon_1 ) H^2 u_1^{(2)}(H)\varepsilon u_1(H)\\&-8(2853 \varepsilon \varepsilon_1+220) H^3 (u_1(H))^3 + 8(12739 \varepsilon\varepsilon_1 + 1435) H^4 u_1^{(2)}(H) u_1(H)\\&-96(51\varepsilon\varepsilon_1+4) u_1^{(3)}(H) H^5+ 972 u_1^{(3)}(H) H (u_1(H))^2- 436 (u_1(H))^3 u_1^{(2)}(H)\\&+16(1476\varepsilon c -144c \varepsilon_1+741\varepsilon \lambda +12\lambda \varepsilon_1) H^5 u_1(H)- 264 u_1^{(4)}(H) H^2 u_1(H)\\&+\varepsilon_1(2010 c+134 \lambda) H (u_1(H))^3 +(24576 \varepsilon \varepsilon_1-20736) H^7 u_1(H)\\& + (52\varepsilon\varepsilon_1+ 580) H u_1(H) (u_1^{(2)}(H))^2+24\varepsilon_1(51c-2 \lambda )u_1^{(3)}(H) H^3.
\end{aligned}
\end{equation*}}
which together with \eqref{5---lemma5.3-proof--6} deduce that
$$
c_1f_1+c_2=0,
$$
where $c_1$ and $c_2$ are polynomials about $H$ and $u_1^{(k)}(H)$, $k=1, 2, \cdots, 5$.

From the expression $[u_i, u_1]=\nabla_{u_i}u_1-\nabla_{u_1}u_i$, as well as the relations \eqref{5---B(H)--4case2} and  $\Gamma_{i1}^1=\Gamma_{1i}^1=0$ (see eq. \eqref{5---compatibility-eq1--4} and \eqref{5---codazzi-eq--1}), we have
\begin{equation}\label{5---symmetry--1}
\begin{aligned}
u_iu_1^{(k)}(H)=0,\ k=0, 1, 2, \cdots, 5.
\end{aligned}
\end{equation}
Differentiate $c_1f_1+c_2=0$ along $u_i$, with $i=2, 3, 4$, using \eqref{5---symmetry--1}, we have
$$
c_1u_i(f_1)=0.
$$

Suppose $u_i(f_1)\neq 0$, for some $i\in\{2, 3, 4\}$, we can conclude that $c_1=c_2=0$. And then, we can eliminate $u_1^{(k)}(H)$, with $k=1, 2, \cdots, 5$ and obtain a polynomial equation about $H$, which implies $H$ is a constant, a contradiction.
So, $u_i(f_1)=0$, $i=2, 3, 4$.
$\hfill\square$

\vskip.2cm
\noindent
{\bf Lemma 3.5}\quad \emph{We have $u_i(\Gamma_{j1}^j)=u_i(\lambda_j)=0$, with $i, j=2, 3, 4$.}

{\bf Proof}\quad Since $u_i(f_1)=0$, with $i=2, 3, 4$ (see Lemma 3.4), it follows
\begin{equation}\label{5---lemma5.4--1}
\begin{aligned}
u_iu_1^{(k)}(f_1)=0,\  i=2, 3, 4,\ k=0, 1, 2, 3.
\end{aligned}
\end{equation}
Differentiate the expressions of $f_k$, $k=2, 3$ in Lemma 3.3, combining
\eqref{5---symmetry--1} and \eqref{5---lemma5.4--1}, we find
\begin{equation*}
\begin{aligned}
u_i(f_k)=0,\  i=2, 3, 4,\ k=1,2,3.
\end{aligned}
\end{equation*}
Note that $f_k=\sum_{i=2}^4(\Gamma_{i1}^i)^k$, the above equation tells us that
\begin{equation*}
\begin{cases}
u_i(\Gamma_{21}^{2})+u_i(\Gamma_{31}^{3})+u_i(\Gamma_{41}^{4})=0,\\
\Gamma_{21}^{2}u_i(\Gamma_{21}^{2})+\Gamma_{31}^{3}u_i(\Gamma_{31}^{3})+\Gamma_{41}^{4}u_i(\Gamma_{41}^{4})=0,\\
(\Gamma_{21}^{2})^2u_i(\Gamma_{21}^{2})+(\Gamma_{31}^{3})^2u_i(\Gamma_{31}^{3})+(\Gamma_{41}^{4})^2u_i(\Gamma_{41}^{4})=0,
\end{cases}
\end{equation*}
with $i=2, 3, 4$. Observe \eqref{5---Gauss-eq--1}, we know $\Gamma_{21}^{2}$, $\Gamma_{31}^{3}$ and $\Gamma_{41}^{4}$ are distinct. So, the coefficient determinant of the above system
\begin{equation*}
\begin{aligned}
\left|
  \begin{array}{ccc}
     1   &1      &1\\
      \Gamma_{21}^{2}   & \Gamma_{31}^{3} & \Gamma_{41}^{4}\\
     (\Gamma_{21}^{2})^2&     (\Gamma_{31}^{3})^2&  (\Gamma_{41}^{4})^2
  \end{array}
\right|=(\Gamma_{41}^{4}-\Gamma_{31}^{3})(\Gamma_{31}^{3}-\Gamma_{21}^{2})(\Gamma_{41}^{4}-\Gamma_{21}^{2})\neq 0.
\end{aligned}
\end{equation*}
So, the above system has only zero solution, i.e. $u_i(\Gamma_{j1}^{j})=0$, $i, j=2, 3, 4$.
And then, we have
$$
u_iu_1(\Gamma_{j1}^{j})=0,\ i, j=2, 3, 4.
$$
Differentiate \eqref{5---Gauss-eq--1} along $u_j$, combining $u_j(\Gamma_{i1}^{i})=0$ and $u_ju_1(\Gamma_{i1}^{i})=0$,
we get $u_j(\lambda_i)=0$, $i, j=2, 3, 4$.
$\hfill\square$

\subsubsection{The proof of Proposition 3.2}

\vskip.2cm
Assume that $H$ is not a constant, then there exists a neighbourhood $U_p$ of $p$  such that $H\neq 0$ and $\nabla H\neq0$.
The equation \eqref{proper-eq1} implies $\nabla H$ is an eigenvector of $A$, with corresponding eigenvalue $-2\varepsilon H$. Without loss of generality, we suppose $\nabla H$ is in the direction of $u_1$ and $\lambda_1=-2\varepsilon H$, then the equations and Lemmas in subsection 3.1.1 hold. We will use these equations and Lemma 3.5 to deduce contradictions.

\vskip.1cm
The second equation in \eqref{5---codazzi-eq--2} yields
\begin{equation}\label{5---codazzi-eq--2-1}
\begin{aligned}
&\varepsilon_4(\lambda_2-\lambda_4)\Gamma_{32}^{4}=\varepsilon_4(\lambda_3-\lambda_4)\Gamma_{23}^{4}=
\varepsilon_3(\lambda_4-\lambda_3)\Gamma_{24}^{3}\\=&\varepsilon_3(\lambda_2-\lambda_3)\Gamma_{42}^{3}
=\varepsilon_2(\lambda_3-\lambda_2)\Gamma_{43}^{2}=
\varepsilon_2(\lambda_4-\lambda_2)\Gamma_{34}^{2}.
\end{aligned}
\end{equation}
Together \eqref{5---compatibility-eq1--4} with \eqref{5---codazzi-eq--2-1}, we have
$$
(\lambda_2-\lambda_3)(\lambda_3-\lambda_4)(\lambda_4-\lambda_2)(\Gamma_{23}^{4}\Gamma_{32}^{4}+\Gamma_{24}^{3}\Gamma_{42}^{3}+\Gamma_{34}^{2}\Gamma_{43}^{2})=0,
$$
i.e.
\begin{equation}\label{5---codazzi-eq--2-2}
\begin{aligned}
\Gamma_{23}^{4}\Gamma_{32}^{4}+\Gamma_{24}^{3}\Gamma_{42}^{3}+\Gamma_{34}^{2}\Gamma_{43}^{2}=0.
\end{aligned}
\end{equation}

It follows from the equation \eqref{5---codazzi-eq--2} and Lemma 3.5 that
\begin{equation}\label{5---gammaiji-1}
\Gamma_{ij}^{i}=0,\ i, j=2, 3, 4.
\end{equation}
Apply Gauss equation for $\langle R(u_i, u_j)u_k, u_1\rangle$ and $\langle R(u_i, u_j)u_i, u_j\rangle$, with $i, j, k$ are distinct and $i, j, k=2, 3, 4$, combining \eqref{5---compatibility-eq1--4}, \eqref{5---codazzi-eq--1} and \eqref{5---gammaiji-1}, we have
\begin{equation}\label{5---Gauss-eq--2}
\begin{aligned}
&\varepsilon_4(\Gamma_{21}^{2}-\Gamma_{41}^{4})\Gamma_{32}^{4}=\varepsilon_4(\Gamma_{31}^{3}-\Gamma_{41}^{4})\Gamma_{23}^{4}=
\varepsilon_3(\Gamma_{41}^{4}-\Gamma_{31}^{3})\Gamma_{24}^{3}\\=&\varepsilon_3(\Gamma_{21}^{2}-\Gamma_{31}^{3})\Gamma_{42}^{3}
=\varepsilon_2(\Gamma_{31}^{3}-\Gamma_{21}^{2})\Gamma_{43}^{2}=\varepsilon_2(\Gamma_{41}^{4}-\Gamma_{21}^{2})\Gamma_{34}^{2},
\end{aligned}
\end{equation}
and
\begin{equation}\label{5---Gauss-eq--3qian}
\begin{cases}
\varepsilon_1\varepsilon_2\varepsilon_3\Gamma_{21}^{2}\Gamma_{31}^{3}-2\varepsilon_4\Gamma_{23}^{4}\Gamma_{32}^{4}=-c\varepsilon_2\varepsilon_3-\lambda_2\lambda_3\varepsilon\varepsilon_2\varepsilon_3,\\
\varepsilon_1\varepsilon_2\varepsilon_4\Gamma_{21}^{2}\Gamma_{41}^{4}-2\varepsilon_3\Gamma_{24}^{3}\Gamma_{42}^{3}=-c\varepsilon_2\varepsilon_4-\lambda_2\lambda_4\varepsilon\varepsilon_2\varepsilon_4,\\
\varepsilon_1\varepsilon_3\varepsilon_4\Gamma_{31}^{3}\Gamma_{41}^{4}-2\varepsilon_2\Gamma_{34}^{2}\Gamma_{43}^{2}=-c\varepsilon_3\varepsilon_4-\lambda_3\lambda_4\varepsilon\varepsilon_3\varepsilon_4,
\end{cases}
\end{equation}
which together with \eqref{5---codazzi-eq--2-2} implies that
\begin{equation}\label{5---Gauss-eq--3}
\varepsilon_1(\Gamma_{21}^{2}\Gamma_{31}^{3}+\Gamma_{21}^{2}\Gamma_{41}^{4}+\Gamma_{31}^{3}\Gamma_{41}^{4})+\varepsilon(\lambda_2\lambda_3+\lambda_3\lambda_4+\lambda_2\lambda_4)
+3c=0.
\end{equation}

In the following, we treat the cases $\Gamma_{23}^{4}\neq 0$ at some point and $\Gamma_{23}^{4}=0$ at any point in $U_p$, respectively.

\vskip.2cm
\emph{Case 1:\ $\Gamma_{23}^{4}\neq 0$ at some point in $U_p$.}

\vskip.1cm
Suppose $\Gamma_{23}^{4}\neq 0$ at $q\in U_p$, then there exists a neighbourhood  $U_q\subset U_p$, such that $\Gamma_{23}^{4}\neq 0$ on $U_q$. We work on $U_q$ in the following discussion.
The equations \eqref{5---codazzi-eq--2-1} and \eqref{5---Gauss-eq--2} give
\begin{equation*}
\frac{\Gamma_{41}^{4}-\Gamma_{31}^{3}}{\lambda_4-\lambda_3}=\frac{\Gamma_{41}^{4}-\Gamma_{21}^{2}}{\lambda_4-\lambda_2}=\frac{\Gamma_{31}^{3}-\Gamma_{21}^{2}}{\lambda_3-\lambda_2}.
\end{equation*}
As $u_i(\Gamma_{j1}^{j})=u_i(\lambda_j)=0$, $i, j=2, 3, 4$, we conclude from the above equation that there exists two smooth functions $\mu$ and $\nu$, with $u_i(\mu)=u_i(\nu)=0$, $i=2, 3, 4$, such that
\begin{equation}\label{5---gamma--4}
\Gamma_{i1}^i=\mu\lambda_i+\nu,\ i=2, 3, 4.
\end{equation}
Take sum for $i$ in \eqref{5---gamma--4}, then
\begin{equation}\label{5---f1--2}
f_1=6\varepsilon H\mu+3\nu.
\end{equation}
Differentiate \eqref{5---gamma--4} along $u_1$, using \eqref{5---codazzi-eq--2} and \eqref{5---Gauss-eq--1},
we get for $i=2, 3, 4$,
$$
(u_1(\mu)-2\varepsilon H\mu^2+\mu\nu-2H\varepsilon_1)\lambda_i+u_1(\nu)+\nu^2-2\varepsilon H\mu\nu+c\varepsilon_1=0,
$$
which yields
\begin{equation}\label{5---u1mu-u1nu-1}
\begin{cases}
u_1(\mu)=2\varepsilon H\mu^2-\mu\nu+2H\varepsilon_1,\\
u_1(\nu)=-\nu^2+2\varepsilon H\mu\nu-c\varepsilon_1.
\end{cases}
\end{equation}
Substitute \eqref{5---gamma--4} into \eqref{5---Gauss-eq--3}, we obtain
\begin{equation*}
(\varepsilon_1\mu^2+\varepsilon)(\lambda_2\lambda_3+\lambda_2\lambda_4+\lambda_3\lambda_4)+2\varepsilon_1\mu\nu(\lambda_2+\lambda_3+\lambda_4)+3\varepsilon_1\nu^2
+3c=0,
\end{equation*}
which can be rewritten as
\begin{equation}\label{5---Gauss-eq--4}
(\varepsilon_1\mu^2+\varepsilon)\text{tr}A^2=
(\varepsilon_1\mu^2+\varepsilon)40H^2+24\varepsilon\varepsilon_1 H\mu\nu+6\varepsilon_1\nu^2
+6c.
\end{equation}
Put \eqref{5---gamma--4} and \eqref{5---f1--2} into \eqref{5---lemma5.2-proof-2}, combining \eqref{5---Gauss-eq--4}, we have
\begin{equation}\label{5---u1H--1}
u_1(H)=\frac{-26}{3}\varepsilon H^2\mu-\frac{4\varepsilon_1H\mu^2\nu+\varepsilon\varepsilon_1\mu\nu^2+c\varepsilon\mu}{\varepsilon_1\mu^2+\varepsilon}
-2H\nu.
\end{equation}

Act on \eqref{5---u1H--1} by $u_1$, applying \eqref{5---u1mu-u1nu-1}, we get the expression of $u_1u_1(H)$. And then, substitute the expressions of $u_1(H)$ and $u_1u_1(H)$ into \eqref{5--proper-td-1}, combining \eqref{5---Gauss-eq--4}, it gives
\begin{equation}\label{5-----1}
\begin{aligned}
&(36 \mu^3\varepsilon + 18\mu\varepsilon_1)\nu ^3 + (300  H \mu^4 + 156 \mu^2 H  \varepsilon\varepsilon_1 - 72H )\nu ^2 + (816 \mu^5 H ^2\varepsilon \\
 &+ 624 \mu^3 H ^2\varepsilon_1+ 36 \mu^3 c\varepsilon\varepsilon_1 - 192\varepsilon H ^2\mu+ 18 c\mu ) \nu+ 940 \mu^4 H ^3\varepsilon\varepsilon_1 \\
  & + 728 \mu^6 H ^3+ 174 \mu^4 H c\varepsilon_1 + 9 \mu^4 H  \lambda\varepsilon_1- 304 \mu^2 H ^3 + 120 \mu^2 H  c\varepsilon\\& + 18 \mu^2 H  \lambda\varepsilon- 516 H ^3\varepsilon\varepsilon_1 - 54 H  c\varepsilon_1 + 9\varepsilon_1 \lambda H=0.
\end{aligned}
\end{equation}
By differentiating \eqref{5-----1} along $u_1$, using \eqref{5---u1mu-u1nu-1} and \eqref{5---u1H--1}, one derive
\begin{equation}\label{5-----2}
\begin{aligned}
d_1\nu^4+d_2\nu^3+d_3\nu^2+d_4\nu+d_5=0,
\end{aligned}
\end{equation}
where
\begin{equation*}
\begin{aligned}
d_1=&-516 \mu ^5 \varepsilon \varepsilon_1 - 444 \mu ^3,\\
d_2=&-4800 \mu ^6 H  \varepsilon_1 - 4416 \mu ^4 H  \varepsilon + 420 \mu ^2 H  \varepsilon_1 + 324 H  \varepsilon,\\
d_3=&-15872 \mu ^7 H ^2 \varepsilon \varepsilon_1 - 17668 \mu ^5 H ^2 - 726 \mu ^5 c \varepsilon - 9 \mu ^5 \lambda \varepsilon + 1864 \mu ^3 H ^2 \varepsilon \varepsilon_1 \\
 &- 618 \mu ^3 c \varepsilon_1- 18 \mu ^3 \lambda \varepsilon_1 + 3660 \mu  H ^2 + 36 c \varepsilon \mu  - 9 \mu  \lambda \varepsilon,\\
d_4=&-21824 \mu ^8 H ^3 \varepsilon_1 - 31432 \mu^6 H ^3 \varepsilon - 3684 \mu ^6 H  c - 90 \mu ^6 H  \lambda + 5320 \mu ^4 H ^3 \varepsilon_1 \\&
- 3588 \mu ^4 H  c \varepsilon \varepsilon_1 - 198 \mu ^4 H  \lambda
\varepsilon \varepsilon_1 + 17640 \mu ^2 H ^3 \varepsilon + 384 \mu ^2 H  c - 126 \mu ^2 H  \lambda\\& + 2712 H ^3 \varepsilon_1 + 288 H  c \varepsilon \varepsilon_1 - 18 H  \lambda \varepsilon \varepsilon_1,\\
d_5=& 6024 \mu ^5 H ^4 \varepsilon \varepsilon_1-10192 \mu ^9 H ^4 \varepsilon \varepsilon_1 - 18376 \mu^7 H ^4 - 3116 \mu ^7 H ^2 c \varepsilon - 6 \mu ^7 H ^2 \lambda \varepsilon\\& - 3544 \mu ^5 H ^2 c \varepsilon_1 - 18 \mu ^5 H ^2 \lambda \varepsilon_1 - 210 \mu ^5 c^2 \varepsilon \varepsilon_1 - 9 \mu ^5 c \lambda \varepsilon \varepsilon_1 + 26408 \mu ^3 H ^4\\& + 2260 \mu ^3 H ^2 c \varepsilon - 18 \mu ^3 H ^2 \lambda \varepsilon + 12200 \mu  H ^4 \varepsilon \varepsilon_1 - 174 \mu ^3 c^2 - 18 \mu ^3 c \lambda \\&+ 2688 \mu  H ^2 c \varepsilon_1 - 6 \mu  H ^2 \lambda \varepsilon_1 + 36 \mu  c^2 \varepsilon \varepsilon_1 - 9 \mu  c \lambda \varepsilon \varepsilon_1.
\end{aligned}
\end{equation*}
Eliminate $\nu$ from \eqref{5-----1} and \eqref{5-----2}, we can get a polynomial equation about $H$ and $\mu$
\begin{equation}\label{5---fHmu}
\begin{aligned}
\sum_{k=29}^{33}\sum_{l=k-28}^7r_{k,l}H^{2l}\mu^{2k}
+\sum_{k=4}^{28}\sum_{l=0}^7r_{k,l}H^{2l}\mu^{2k}+\sum_{k=0}^{3}\sum_{l=4-k}^7r_{k,l}H^{2l}\mu^{2k}=0,
\end{aligned}
\end{equation}
where $r_{0,l}, r_{1,l}, \cdots, r_{33,l}$ are all real constants.

Applying $u_1$ to both sides of \eqref{5-----2}, we obtain a polynomial equation about $H, \mu$ and $\nu$.
Then, eliminating $\nu$ from this equation and \eqref{5-----1}, we derive another polynomial equation about $H$ and $\mu$
\begin{equation*}
\begin{aligned}
(\varepsilon\varepsilon_1 + \mu^2)^{10} (\varepsilon\varepsilon_1 + 2 \mu^2)^6(g_1(H,\mu))^2g_2(H,\mu)=0,
\end{aligned}
\end{equation*}
i.e.
\begin{equation}\label{5---gHmu}
\begin{aligned}
(\varepsilon\varepsilon_1 + \mu^2) (\varepsilon\varepsilon_1 + 2 \mu^2)g_1(H,\mu)g_2(H,\mu)=0,
\end{aligned}
\end{equation}
where
$$
\begin{aligned}
g_1(H,\mu)= &6160\varepsilon_1 H^2 \mu^{14} +
 16(17136 c - 62986 \varepsilon H^2  +
1935 \lambda) \mu^{12}+69984\varepsilon H^2\\
 & + 8\varepsilon_1(37935 \varepsilon c
 - 512645 H^2 + 8424 \varepsilon
 \lambda) \mu^{10}+
 48(987\lambda  - 28138\varepsilon H^2\\
 & +
 1884 c) \mu^8+ \varepsilon_1(8334 \varepsilon c
 + 2136896  H^2  + 10044 \varepsilon
 \lambda) \mu^6 +
 ( 6570 c\\
 & +312180 \varepsilon H^2-
 1467 \lambda)\mu^4+ \varepsilon_1 (2106 \varepsilon c
 - 14292 H^2  -567 \varepsilon
 \lambda) \mu^2,
 \end{aligned}
$$
and
$$
\begin{aligned}
g_2(H,\mu)\!=\!&\sum_{l=4}^6s_{17,l}H^{2l}\mu^{34}+\sum_{l=3}^6s_{16,l}H^{2l}\mu^{32}+\sum_{l=2}^6s_{15,l}H^{2l}\mu^{30}+\sum_{l=1}^6s_{14,l}H^{2l}\mu^{28}\\
&+\sum_{k=2}^{13}\sum_{l=0}^6s_{k,l}H^{2l}\mu^{2k}+\sum_{l=1}^6s_{1,l}H^{2l}\mu^{2}+\sum_{l=2}^6s_{0,l}H^{2l},
\end{aligned}
$$
with $s_{0,l}, s_{1,l}, \cdots, s_{17,l}$ are all real constants.
We can eliminate $\mu$ from \eqref{5---fHmu} and $\eqref{5---gHmu}$, and give a polynomial equation of degree $696$ for $H$, which implies $H$ is a constant in $U_q$. Thus, $\nabla H=0$ at $q$, a contradiction.

\vskip.2cm
\emph{Case 2:\ $\Gamma_{23}^{4}\equiv0$ on $U_p$.}

\vskip.1cm
In this case, \eqref{5---Gauss-eq--3qian} can be simplified into
\begin{equation}\label{5-II--Gauss-eq--3qian}
\begin{cases}
\varepsilon_1\Gamma_{21}^{2}\Gamma_{31}^{3}=-c-\varepsilon\lambda_2\lambda_3,\\
\varepsilon_1\Gamma_{21}^{2}\Gamma_{41}^{4}=-c-\varepsilon\lambda_2\lambda_4,\\
\varepsilon_1\Gamma_{31}^{3}\Gamma_{41}^{4}=-c-\varepsilon\lambda_3\lambda_4.
\end{cases}
\end{equation}
As $\lambda_2, \lambda_3$ and $\lambda_4$ are distinct, we have $c+\varepsilon\lambda_2\lambda_3\neq 0$, or $c+\varepsilon\lambda_2\lambda_4\neq0$, or $c+\varepsilon\lambda_3\lambda_4\neq 0$.
So, the above equations implies there are at least two non-zero terms in $\{\Gamma_{21}^{2}, \Gamma_{31}^{3}, \Gamma_{41}^{4}\}$. Without loss of generality, we suppose $\Gamma_{21}^{2}, \Gamma_{31}^{3}\neq 0$ at any point in $U_p$, and discuss separately the subcases $\Gamma_{41}^{4}\equiv0$ on $U_p$ and $\Gamma_{41}^{4}\neq 0$ at some point.

\vskip.2cm
\emph{Subcase (1):\ $\Gamma_{41}^{4}\equiv0$ on $U_p$.}

\vskip.1cm
It follows from \eqref{5-II--Gauss-eq--3qian} that $\lambda_4=0$.
If $c\neq 0$, then $c+\lambda_2\lambda_4\neq 0$, which implies $\Gamma_{41}^{4}\neq 0$, a contradiction.
In the following, we suppose $c=0$.
Act on both sides of $\lambda_2+\lambda_3=6\varepsilon H$ by $u_1$, we get
\begin{equation}\label{5---II.2--u1H}
6\varepsilon u_1(H)=(-2\varepsilon H-\lambda_2)\Gamma_{21}^2+(-2\varepsilon H-\lambda_3)\Gamma_{31}^3.
\end{equation}
Differentiate \eqref{5---II.2--u1H} along $u_1$, using \eqref{5---codazzi-eq--2} and \eqref{5---Gauss-eq--1}, we obtain the expression of $u_1u_1(H)$, and put it into \eqref{5--proper-td-1}, combining \eqref{5-II--Gauss-eq--3qian} and \eqref{5---II.2--u1H}, we have
\begin{equation}\label{5---II.2--propertdx2}
\begin{aligned}
2(-2\varepsilon H-\lambda_2)(\Gamma_{21}^2)^2&+2(-2\varepsilon H-\lambda_3)(\Gamma_{31}^3)^2-34\varepsilon_1H\lambda_2\lambda_3\\&
+468\varepsilon_1H^3-9\varepsilon\varepsilon_1\lambda H=0.
\end{aligned}
\end{equation}
By applying \eqref{5---codazzi-eq--2}, \eqref{5---Gauss-eq--1} and \eqref{5---II.2--u1H}, differentiate  \eqref{5---II.2--propertdx2} along $u_1$,
and combining \eqref{5-II--Gauss-eq--3qian} and \eqref{5---II.2--propertdx2}, we derive
\begin{equation}\label{5---II.2--upropertdx2}
\begin{aligned}
K_1\Gamma_{21}^{2}+K_2\Gamma_{31}^3=0,
\end{aligned}
\end{equation}
where
\begin{equation*}
\begin{aligned}
K_1=&-680 H  \lambda_2  \lambda_3  \varepsilon \varepsilon_1 + 6552 \varepsilon H ^3 - 96 H ^2 \lambda_2  \varepsilon + 264 H  \lambda_2
\lambda_3  \varepsilon \\&+ 40 \lambda_2 ^2 \lambda_3  \varepsilon - 1404 H ^2 \lambda_2
  + 408 H ^2 \lambda_3  - 48 H  \lambda_2 ^2 + 80 H  \lambda_2  \lambda_3  \\&+ 34 \lambda_3  \lambda_2 ^2 - 4 \lambda_2  \lambda_3 ^2+ 9 \lambda_2  \lambda \varepsilon - 162 H  \lambda;\\
K_2=&-680 H  \lambda_2  \lambda_3  \varepsilon \varepsilon_1 + 6552 \varepsilon H ^3 - 96 H ^2 \lambda_3  \varepsilon + 264 H  \lambda_2
\lambda_3  \varepsilon\\& + 40 \lambda_3 ^2 \lambda_2  \varepsilon - 1404 H ^2 \lambda_3
 + 408 H ^2 \lambda_2  - 48 H  \lambda_3 ^2 + 80 H  \lambda_2  \lambda_3 \\& + 34 \lambda_2  \lambda_3 ^2 - 4 \lambda_3  \lambda_2 ^2+ 9 \lambda_3  \lambda \varepsilon - 162 H  \lambda.
\end{aligned}
\end{equation*}
Since \eqref{5-II--Gauss-eq--3qian}, we can eliminate $\Gamma_{21}^{2}$ and $\Gamma_{31}^3$ from \eqref{5---II.2--propertdx2} and \eqref{5---II.2--upropertdx2}, and give
\begin{equation}\label{5---II.2--lambda23H---1}
\begin{aligned}
&-32(\lambda_2^6\lambda_3^3+\lambda_2^3\lambda_3^6)+c_1(\lambda_2^5\lambda_3^4+\lambda_2^4\lambda_3^5)+c_2(\lambda_2^5\lambda_3^3+\lambda_2^3\lambda_3^5)
+c_3\lambda_2^4\lambda_3^4\\&+c_4(\lambda_2^5\lambda_3^2+\lambda_2^2\lambda_3^5)
+c_5(\lambda_2^4\lambda_3^3+\lambda_2^3\lambda_3^4)
-9216\varepsilon H^3(\lambda_2^5\lambda_3+\lambda_2\lambda_3^5)\\&+c_6(\lambda_2^4\lambda_3^2+\lambda_2^2\lambda_3^4)
+c_7\lambda_2^3\lambda_3^3+c_8(\lambda_2^4\lambda_3+\lambda_2\lambda_3^4)+c_{9}(\lambda_2^3\lambda_3^2+\lambda_2^2\lambda_3^3)\\&
+c_{10}(\lambda_2^3\lambda_3+\lambda_2\lambda_3^3)+c_{11}\lambda_2^2\lambda_3^2+c_{12}(\lambda_2^3+\lambda_3^3)
+c_{13}(\lambda_2^2\lambda_3+\lambda_2\lambda_3^2)\\&
+c_{14}(\lambda_2^2+\lambda_3^2)+c_{15}\lambda_2\lambda_3+c_{16}(\lambda_2+\lambda_3)+c_{17}=0,
\end{aligned}
\end{equation}
where $c_i$, $1\leq i\leq 17$ are polynomials about $H$.
Act on \eqref{5---II.2--lambda23H---1} by $u_1$, we get
\begin{equation}\label{5---II.2--upropertdx2--l}
\begin{aligned}
L_1\Gamma_{21}^{2}+L_2\Gamma_{31}^3=0,
\end{aligned}
\end{equation}
where $L_1$ and $L_2$ are polynomials about $H, \lambda_2$ and $\lambda_3$.
It follows from \eqref{5---II.2--upropertdx2} and \eqref{5---II.2--upropertdx2--l} that
$$
K_1L_2-K_2L_1=0.
$$
Take into account $\lambda_2+\lambda_3=6\varepsilon H$, we can eliminate $\lambda_2$ and $\lambda_3$ from \eqref{5---II.2--lambda23H---1} and the above equation, and derive a polynomial equation of degree $78$ for $H$, which implies $H$ is a constant, a contradiction.

\vskip.2cm
\emph{Subcase (2):\ $\Gamma_{41}^{4}\neq 0$ at some point in $U_p$.}

\vskip.1cm
Suppose $\Gamma_{41}^{4}\neq 0$ at $q\in U_p$, then there exists a neighbourhood  $U_q\subset U_p$ such that $\Gamma_{41}^{4}\neq 0$ on $U_q$. For this subcase, we work on $U_q$.
Let
$$\mu=\lambda_2\lambda_3+\lambda_2\lambda_4+\lambda_3\lambda_4$$
and
$$
\nu=\lambda_2\lambda_3\lambda_4.
$$
By use of \eqref{5---codazzi-eq--2}, differentiate both sides of $\lambda_2+\lambda_3+\lambda_4=6\varepsilon H$, $\mu=\lambda_2\lambda_3+\lambda_2\lambda_4+\lambda_3\lambda_4$ and $
\nu=\lambda_2\lambda_3\lambda_4$ along $u_1$, and then multiply these results by $(c+\varepsilon\lambda_3\lambda_4)\Gamma_{21}^2$, combining \eqref{5-II--Gauss-eq--3qian},
we obtain
\begin{equation}\label{5-II--u1H-u1mu-u1nu}
\begin{cases}
6\varepsilon u_1(H)=Q(-12\varepsilon Hc^2-10H\mu c+3\varepsilon\nu c-48H^2\nu+2\mu\nu),\\ u_1(\mu)=Q[3\nu^2- 2(12H^2+\mu)c^2-2(\mu^2+6\varepsilon H^2\mu-3H\nu)c-10\varepsilon H\mu\nu],\\ u_1(\nu)=Q[(-2\varepsilon H\mu-3\nu)c^2+(-24\varepsilon H^2\nu-2\varepsilon\mu\nu)c-12\varepsilon H\nu^2],
\end{cases}
\end{equation}
where $Q=\frac{-\varepsilon_1}{(c+\varepsilon\lambda_3\lambda_4)\Gamma_{21}^2}$.
As \eqref{5-II--Gauss-eq--3qian} and \eqref{5-II--u1H-u1mu-u1nu},
\eqref{5--proper-td-1} can be rewritten as
\begin{equation}\label{5--proper-td-2}
\begin{aligned}
&(30 \varepsilon H  c^2 - 44 H  c^2 - 6 \varepsilon \nu  c - 2 c \nu ) \mu ^2 + (504 H ^3 c^2 - 36 H  c^3 \varepsilon\\& - 9 H  c^2 \lambda \varepsilon - 12 c \nu  H ^2- 15 \nu c^2 \varepsilon - 78 H  \nu ^2) \mu+ 3024 H ^4 \nu  c\\& + 504 H ^3 \varepsilon c^3- 108 H ^2 \nu  \varepsilon c^2 - 54 H ^2 \nu  \lambda \varepsilon c+ 1080 H ^3 \nu
 ^2 \\&- 72 H  \nu ^2 \varepsilon c- 9 H  \nu ^2 \lambda \varepsilon - 18 H  c^4 - 9 H \lambda c^3 + 9 \nu ^3 \varepsilon - 9 \nu  c^3=0.
\end{aligned}
\end{equation}
Act on \eqref{5--proper-td-2} by $u_1$, applying \eqref{5-II--u1H-u1mu-u1nu},
it gives
\begin{equation}\label{5--u1proper-td-2}
\begin{aligned}
f_1\mu^3+f_2\mu^2 +f_3\mu +f_4=0.
\end{aligned}
\end{equation}
where
\begin{equation*}
\begin{aligned}
f_1=&(828 \varepsilon  -256) H c^3 +(180\varepsilon + 80)\nu c^2;\\
f_2=& 4(522\varepsilon + 747)H\nu^2 c-1728(6\varepsilon + 5)H^3 c^3 - 288(11\varepsilon - 46)H^2\nu c^2 - 156\nu^3\\&+ 4(531 \varepsilon - 162)H c^4 +18 (5\varepsilon+ 6) H\lambda c^3 + 6(10\varepsilon+ 87)\nu c^3- 18\nu c^2\lambda\varepsilon ;\\
f_3=&-72576H^5 c^3 - (22752\varepsilon+ 8640)H^3 c^4 + 1296H^3 c^3\lambda\varepsilon - 283392H^4\nu c^2
\\& + 2916H^2\nu c^2\lambda\varepsilon- 3312H^3\nu^2 c + (17928\varepsilon+ 3888)H^2\nu c^3 + 90H\nu^2 c\lambda\varepsilon\\& + (6552\varepsilon+ 648)H\nu^2 c^2 + 1152H c^5
 + 306H c^4\lambda + 26136H^2\nu^3 \\&- 774\nu^3 c\varepsilon - 18\nu^3\lambda\varepsilon + 414\nu c^4 - 45\nu c^3\lambda - 216\nu^3 c;
 \end{aligned}
 \end{equation*}
\begin{equation*}
\begin{aligned}
 f_4=&7776H^4\nu\lambda\varepsilon c^2 + 11664H^3\nu^2\lambda\varepsilon c-311040H^4\nu^3 - 236736H^4\nu\varepsilon c^3\\& + 17280H^3\nu^2\varepsilon c^2 + 108H\lambda\varepsilon c^5 + 16632H^2\nu^3\varepsilon c + 1728H^2\nu^3\lambda\varepsilon - 27\nu\lambda\varepsilon c^4\\& - 54H\nu^2\lambda c^2 + 2376H^2\nu\lambda c^3 - 1109376H^5\nu^2 c + 12096H^2\nu c^4- 27\nu^3\lambda c\\& - 435456H^6\nu c^2 + 1296H^3\lambda c^4 + 2268H\nu^2 c^3- 72576H^5\varepsilon c^4 + 216H\varepsilon c^6\\& + 108\nu\varepsilon c^5 - 3348H\nu^4\varepsilon - 12960H^3 c^5 - 972\nu^3 c^2.
\end{aligned}
\end{equation*}
By differentiating \eqref{5--u1proper-td-2} along $u_1$, using \eqref{5-II--u1H-u1mu-u1nu},
we have another polynomial equation about $H$, $\mu$ and $\nu$, denoted by $g(H, \mu, \nu)=0$.
Eliminating $\mu$ from \eqref{5--proper-td-2}, \eqref{5--u1proper-td-2} and $g(H, \mu, \nu)=0$, we get
\begin{equation}\label{5---eq1-Hnu}
\begin{aligned}
(24336\varepsilon + 18928)\nu ^{15}+\sum_{i=0}^{14}h_i\nu ^i=0,
\end{aligned}
\end{equation}
and
\begin{equation}\label{5---eq2-Hnu}
\begin{aligned}
32(32436060H^2\varepsilon +50162468H^2 +2184345c\varepsilon +1894023 c)\nu^{19}+\sum_{i=0}^{18}l_i\nu ^i=0,
\end{aligned}
\end{equation}
where $h_i$, with $0\leq i\leq 14$ and $l_j$, with $0\leq j\leq 18$ are polynomials about $H$.

When $\varepsilon=1$, \eqref{5---eq1-Hnu} and \eqref{5---eq2-Hnu} can be rewritten as
$$
\begin{cases}
(2cH-\nu)(7cH+4\nu)^2(c^3-6cH\nu-2\nu^2)^3h_1(H,\nu)=0,\\
(2cH-\nu)(7cH+4\nu)^3(c^3-6cH\nu-2\nu^2)^3h_2(H,\nu)=0,
\end{cases}
$$
where
$$
\begin{cases}
h_1(H,\nu)=338\nu^6+\sum_{i=0}^{5}t_i\nu^i,\\
h_2(H,\nu)=(254898c+5162408H^2)\nu^9+\sum_{i=0}^{8}s_i\nu^i.
\end{cases}
$$
Here, $t_i$ and $t_j$ are polynomials about $H$, with $0\leq i\leq 5$ and $0\leq j\leq 8$.
If
$$
(2cH-\nu)(7cH+4\nu)(c^3-6cH\nu-2\nu^2)\neq 0
$$
at some point $o$ in $U_q$,
we can eliminate $\nu$ from $h_1(H,\nu)=0$ and $h_2(H,\nu)=0$, and get a polynomial equation of degree $94$ for $H$, which implies $\nabla H=0$ at $o$, a contradiction.

If
$$
(2cH-\nu)(7cH+4\nu)(c^3-6cH\nu-2\nu^2)\equiv0
$$
on $U_q$, then acting on it by $u_1$, by use of \eqref{5-II--u1H-u1mu-u1nu},
we have
\begin{equation}\label{bu1--sec3.1}
\begin{aligned}
&54 c^6 H \nu -336 c^7 H^2 + 3048 c^5 H^3 \nu +
 12096 c^4 H^5 \nu + 2316 c^4 H^2 \nu^2\\&+ 147 c^5 \nu^2 +
 27936 c^3 H^4 \nu^2 - 1080 c^3 H \nu^3 - 1344 c^2 H^3 \nu^3-
 510 c^2 \nu^4 \\
 &- 10416 c H^2 \nu^4-
 2304 H \nu^5+ ( 98 c^4 \nu^2-292 c^6 H^2 +1008 c^5 H^4\\&+
    130 c^5 H \nu + 4344 c^4 H^3 \nu +
    200 c^3 H^2 \nu^2 - 1532 c^2 H \nu^3 - 340 c \nu^4)\mu =0.
\end{aligned}
\end{equation}
Furthermore, we claim that 
$7cH+4\nu\equiv0$ on $U_q$. Because if $7cH+4\nu\neq0$ (i.e. $(2cH-\nu)(c^3-6cH\nu-2\nu^2)=0$) at some point in $U_q$, then as \eqref{bu1--sec3.1}, we have
\begin{equation}\label{bueq2---sec5}
(c+\lambda_2\lambda_3)(c+\lambda_2\lambda_4)(c+\lambda_3\lambda_4)=\nu^2+6cH\nu+c^2\mu+c^3=0,
\end{equation}
a contradiction.
We can easily derive a polynomial equation about $H$ from $7cH+4\nu=0$ and \eqref{bu1--sec3.1}. Then, $H$ is a constant, a contradiction.

When $\varepsilon=-1$, \eqref{5---eq1-Hnu} and \eqref{5---eq2-Hnu} reduces to
$$
\begin{cases}
(37cH-2\nu)^2L_1(H,\nu)=0,\\
(37cH-2\nu)^3L_2(H,\nu)=0,
\end{cases}
$$
where
$$
\begin{cases}
L_1(H,\nu)=-1352\nu^{13}+\sum_{i=0}^{12}r_i\nu^i,\\
L_2(H,\nu)=(254898c+5162408H^2)\nu^{16}+\sum_{i=0}^{15}w_i\nu^i.
\end{cases}
$$
Here, $r_i$ and $w_j$ are all polynomials about $H$, with $0\leq i\leq 12$ and $0\leq j\leq 8$.
If $37cH-2\nu\neq 0$ at some point $o$ in $U_q$,
we can eliminate $\nu$ from $L_1(H,\nu)=0$ and $L_2(H,\nu)=0$, and get a polynomial equation of degree $292$ for $H$, which implies $\nabla H=0$ at $o$, a contradiction. If $37cH-2\nu=0$ at any point in $U_q$, then acting on it by $u_1$, combining \eqref{5-II--u1H-u1mu-u1nu},
we have
$$
-1488 c \nu H^2 - 346 H \mu c^2 - 444 H c^3 + 144 H\nu^2 +
 98\nu \mu c + 147\nu c^2=0,
$$
which together with $37cH-2\nu=0$ and \eqref{5--proper-td-2} implies that $H$ is a constant in $U_q$, a contradiction.
$\hfill\square$

\subsection{The shape operator has form (I\!I)}

\vskip.2cm
\noindent
{\bf Proposition 3.6}\quad \emph{Let $M^4_r$ be a nondegenerate hypersurface of $N^{5}_s(c)$ with proper mean curvature vector field. Suppose that the shape operator $A$ of $M^4_r$ has the form (I\!I),
then $M^4_r$ has constant mean curvature.}

{\bf Proof}\quad
We suppose $\lambda_1, \lambda_2$ and $\lambda_3$ are distinct with each other.
When there are at most two distinct values in $\{\lambda_1, \lambda_2, \lambda_3\}$,
 we omit the proof, since it is similar but much easier.

Denote $u_2=u_{2_1}$ and $u_3=u_{3_1}$.
Let $\nabla_{u_B}u_C=\Gamma^D_{BC}u_D$, $B, C=1_1, 1_2, 2, 3$,
then compatibility condition implies that
\begin{equation}\label{4.1---compatibility-eq1--4}
\Gamma_{D1_1}^{1_2}=\Gamma_{D1_2}^{1_1}=\Gamma_{D2}^{2}=\Gamma_{D3}^{3}=0,
\end{equation}
and
\begin{equation}\label{4.1---compatibility-eq2--4}
\begin{aligned}
\Gamma_{D1_1}^{1_1}=-\Gamma_{D1_2}^{1_2}, \Gamma_{D2}^{3}=-\varepsilon_2\varepsilon_3\Gamma_{D3}^{2},
\Gamma_{Di}^{1_1}=-\varepsilon_1\varepsilon_i\Gamma_{D1_{2}}^{i},
\Gamma_{Di}^{1_2}=-\varepsilon_1\varepsilon_i\Gamma_{D1_{1}}^{i},
\end{aligned}
\end{equation}
with $i=2, 3$ and $D=1_1, 1_2, 2, 3$. Express $\nabla H$ as
\begin{equation}\label{4.1---nablaH}
\nabla H=\varepsilon_1u_{1_2}(H)u_{1_1}+\varepsilon_1u_{1_1}(H)u_{1_2}+\varepsilon_2u_2(H)u_2+\varepsilon_3u_3(H)u_3.
\end{equation}

Assume that $H$ is not a constant, it follows from \eqref{proper-eq1} that $\nabla H$ is an eigenvector of $A$ with corresponding eigenvalue $-2\varepsilon H$, is light-like or not.

\vskip.2cm
\emph{Case 1:\ $\nabla H$ is light-like.}

\vskip.1cm
In this case, $\nabla H$ is in the direction $u_{1_2}$ and $\lambda_1=-2\varepsilon H$.
As \eqref{4.1---nablaH}, we know
\begin{equation}\label{4.1-case1--B(H)--4case2}
u_{1_{1}}(H)\neq0,\ u_{1_2}(H)=u_2(H)=u_3(H)=0.
\end{equation}
According to symmetry of the connection $\nabla$, we have
\begin{equation}\label{4.1-case1--symmetry--4case2}
\Gamma^{1_{1}}_{BC}=\Gamma^{1_{1}}_{CB},\quad B, C= 1_2, 2, 3.
\end{equation}
Investigate the equation \eqref{Codazzi-td-eq},
for $(X, Y, Z)=(u_{1_1}, u_{i}, u_{1_2}), (u_2, u_3, u_{1_2}),\\ (u_{1_2}, u_{i}, u_{i})$, with $i=2, 3$, by use of \eqref{4.1-case1--B(H)--4case2} and \eqref{4.1-case1--symmetry--4case2}, we obtain that
       \begin{equation}\label{4.1-case1--codazzitds1}
\Gamma_{1_12}^{1_1}=\Gamma_{1_13}^{1_1}=\Gamma_{23}^{1_1}=0,
      \end{equation}
and
\begin{equation}\label{4.1-case1--codazzitds2}
u_{1_2}(\lambda_i)=(-2\varepsilon H-\lambda_i)\Gamma_{i1_2}^{i},\ i=2, 3.
\end{equation}
Calculate $\langle R(u_{1_{2}}, u_i)u_{1_2},u_i\rangle$, $i=2, 3$ by Gauss equation,
combining \eqref{4.1---compatibility-eq1--4}, \eqref{4.1---compatibility-eq2--4}, \eqref{4.1-case1--symmetry--4case2} and \eqref{4.1-case1--codazzitds1}, we get
\begin{equation}\label{4.1-case1--Gausstds1}
u_{1_2}(\Gamma_{i1_2}^{i})
=\Gamma^{1_2}_{1_21_2}
\Gamma_{i1_2}^{i}
-(\Gamma_{i 1_2}^{i})^2,\ i=2, 3.
\end{equation}
As $\lambda_1=-2\varepsilon H$, the equation $2\lambda_1+\lambda_2+\lambda_3=4\varepsilon H$ (cf. eq. \eqref{trA--td---sec2.3}) reduces to
$
\lambda_2+\lambda_3=8\varepsilon H.
$
Differentiate this equation along $u_{1_2}$, applying \eqref{4.1-case1--B(H)--4case2} and \eqref{4.1-case1--codazzitds2},
we find
\begin{equation*}
(-2\varepsilon H-\lambda_2)\Gamma_{21_2}^{2}+(-2\varepsilon H-\lambda_3)\Gamma_{31_2}^{3}=0.
\end{equation*}
By use of \eqref{4.1-case1--B(H)--4case2}, \eqref{4.1-case1--codazzitds2} and \eqref{4.1-case1--Gausstds1}, act $u_{1_2}$ on the above equation,
we have
\begin{equation*}
(-2\varepsilon H-\lambda_2)(\Gamma_{21_2}^{2})^2+(-2\varepsilon H-\lambda_3)(\Gamma_{31_2}^{3})^2=0.
\end{equation*}
Because of $\lambda_2, \lambda_3\neq -2\varepsilon H$ and $\lambda_2+\lambda_3=8\varepsilon H$, we conclude from the above two equations that
\begin{equation}\label{4.1-case1--diff12tds}
\Gamma_{21_2}^{2}=\Gamma_{31_2}^{3}
\end{equation}
and
$$
-12\varepsilon H\Gamma_{21_2}^{2}=0,
$$
which together with $H\neq0$ implies $\Gamma_{21_2}^{2}=0$. Combining \eqref{4.1---compatibility-eq1--4}, \eqref{4.1---compatibility-eq2--4}, \eqref{4.1-case1--symmetry--4case2}, \eqref{4.1-case1--codazzitds1} and \eqref{4.1-case1--diff12tds}, compute $\langle R(u_{1_1}, u_i)u_{1_2}, u_i\rangle$
by using Gauss equation, it gives
$$
2H\lambda_i=c,\ i=2, 3,
$$
which yields that
$
\lambda_2=\lambda_3
$,
a contradiction.

\vskip.2cm
\emph{Case 2:\ $\nabla H$ is not light-like.}

\vskip.1cm
Observe the form (I\!I), we know $\nabla H$ is in the direction $u_2$ or $u_3$ in this case.
Without loss of generality, we suppose $\nabla H$ is in the direction $u_2$, then $\lambda_2=-2\varepsilon H$,
which together with $2\lambda_1+\lambda_2+\lambda_3=4\varepsilon H$ tells us that
$\lambda_3=6\varepsilon H-2\lambda_1$.
The equation \eqref{4.1---nablaH} and the symmetry of connection $\nabla$ give that
\begin{equation}\label{4.1-case2--B(H)--4case2}
u_2(H)\neq0,\ u_B(H)=0,\ B=1_1, 1_2, 3.
\end{equation}
and
\begin{equation}\label{4.1-case2--symmetry--4case2}
\Gamma^2_{BC}=\Gamma^2_{CB},\quad B, C=1_1, 1_2, 3.
\end{equation}
By use of \eqref{4.1---compatibility-eq1--4}, \eqref{4.1---compatibility-eq2--4}, \eqref{4.1-case2--B(H)--4case2}  and \eqref{4.1-case2--symmetry--4case2}, we deduce from \eqref{Codazzi-td-eq} that
\begin{equation}\label{4.1-case2--1--3.1case1}
u_{1_2}(\lambda_1)=\Gamma_{1_23}^{1_1}=\Gamma_{1_23}^{1_1}=\Gamma_{31_1}^{2}=\Gamma_{31_2}^{2}=\Gamma_{31_2}^{3}=\Gamma_{2B}^{2}=0,
\end{equation}
\begin{equation}\label{4.1-case2--1b--3.1case12b--3.1case1}
u_{1_2}(\lambda_3)=(\lambda_1-\lambda_3)\Gamma_{31_2}^{3}, \ \
\Gamma_{1_13}^{1_1}=\Gamma_{1_23}^{1_2}=\frac{u_{3}(\lambda_1)}{\lambda_3-\lambda_1},
\end{equation}
and
\begin{equation}\label{4.1-case2--3--3.1case1}
\Gamma_{1_12}^{1_1}=\Gamma_{1_22}^{1_2}=\frac{u_{2}(\lambda_1)}{-2\varepsilon H-\lambda_1},\ \
\Gamma_{32}^{3}=\frac{u_{2}(\lambda_3)}{-2\varepsilon H-\lambda_3}.
\end{equation}
Let $e_1=\frac{u_{1_1}-u_{1_2}}{\sqrt{2}}$ and $e_2=\frac{u_{1_1}+u_{1_2}}{\sqrt{2}}$,
then $\mathfrak{E}=\{e_1, e_2, u_2, u_3\}$ is an orthonormal basis of $T_xM^4_r$.
We easily find
$$
\nabla_{e_2}e_2(H)-\nabla_{e_1}e_1(H)=\Gamma_{1_11_2}^2+\Gamma_{1_21_1}^2,
$$
and the equation (\ref{proper-eq2}) can be reduced to
\begin{equation}\label{4.1-case2--proper-td--3.1}
u_2u_2(H)+(2\Gamma_{1_12}^{1_1}+\Gamma_{32}^{3})u_2(H)-\varepsilon \varepsilon_2H\text{tr}A^2+\varepsilon_2\lambda H=0,
\end{equation}
where
\begin{equation*}
\text{tr}A^2=2\lambda_1^2+\lambda_3^2+4H^2.
\end{equation*}
Using Gauss equation for
$\langle R(u_2,u_3)u_3,u_2\rangle$ and $\langle R(u_2,u_{1_1})u_{1_2},u_2\rangle$, combining
(\ref{4.1---compatibility-eq1--4}), (\ref{4.1---compatibility-eq2--4}),
(\ref{4.1-case2--symmetry--4case2}) and \eqref{4.1-case2--1--3.1case1}, it is not difficult to check
\begin{equation}\label{4.1-case2-i-j--3.1case1}
\begin{cases}
u_2(\Gamma_{32}^{3})+(\Gamma_{32}^{3})^2=2H\varepsilon_2 \lambda_3-c\varepsilon_2,\\
u_2(\Gamma_{1_12}^{1_1})+(\Gamma_{1_12}^{1_1})^2=2H\varepsilon_2\lambda_1-c\varepsilon_2.
\end{cases}
\end{equation}

Let
$$
f_k=(\Gamma_{1_12}^{1_1})^k+(\Gamma_{1_22}^{1_2})^k+(\Gamma_{32}^{3})^k,\ k=1, 2, \cdots,
$$
then \eqref{5---fkrelations} also holds. We can write $f_k$, $k=2, 3, 4, 5$, as the expressions about $f_1$, $H$, and their derivatives along $u_2$, similarly with \eqref{5---lemma5.2--1}, just need to replace $u_1$ with $u_2$ in \eqref{5---lemma5.2--1}. And then, follow the process of the proof for Lemma 3.4, we have $u_{B}(f_1)=0$, with $B=1_1, 1_2, 3$. Furthermore, we easily get that $u_Bu_2^{(k)}(f_1)=u_Bu_2^{(k)}(H)=0$, with $B=1_1, 1_2, 3$ and $k=1, 2$. Thus, we derive
\begin{equation*}
\begin{aligned}
u_B(f_k)=0,\ B=1_1, 1_2, 3,\ \ k=1,2,
\end{aligned}
\end{equation*}
which gives
\begin{equation*}
\begin{cases}
2u_B(\Gamma_{1_12}^{1_1})+u_B(\Gamma_{32}^{3})=0,\\
2\Gamma_{1_12}^{1_1}u_B(\Gamma_{1_12}^{1_1})+\Gamma_{32}^{3}u_B(\Gamma_{32}^{3})=0.
\end{cases}
\end{equation*}
Observe \eqref{4.1-case2-i-j--3.1case1}, we know $\Gamma_{1_12}^{1_1}\neq\Gamma_{32}^{3}$. So, the above equations yields
$$
u_B(\Gamma_{1_12}^{1_1})=u_B(\Gamma_{32}^{3})=0,\ B=1_1, 1_2, 3.
$$
Act on the equations in \eqref{4.1-case2-i-j--3.1case1} by $u_B$, $B=1_1, 1_2, 3$, using the above equation,
we find
$$
u_B(\lambda_1)=u_B(\lambda_3)=0,\ B=1_1, 1_2, 3,
$$
which together with \eqref{4.1-case2--1b--3.1case12b--3.1case1} implies
\begin{equation}\label{uBlambda10}
\Gamma_{1_13}^{1_1}=\Gamma_{31_2}^{3}=0.
\end{equation}

Applying Gauss equation for
$\langle R(u_{3},u_{1_1})u_3, u_{1_2}\rangle$, combining
(\ref{4.1---compatibility-eq1--4}), (\ref{4.1---compatibility-eq2--4}),
(\ref{4.1-case2--symmetry--4case2}), \eqref{4.1-case2--1--3.1case1} and \eqref{uBlambda10}, we have
\begin{equation}\label{4.1-case2--k--3.1case1}
\Gamma_{32}^{3}\Gamma_{1_12}^{1_1}=2\varepsilon\varepsilon_2\lambda_1^2-6\varepsilon_2H\lambda_1-c\varepsilon_2.
\end{equation}
Substituting $\Gamma_{32}^{3}=\frac{u_2(3\varepsilon H-\lambda_1)}{-4\varepsilon H+\lambda_1}$
and $\Gamma_{1_12}^{1_1}=\frac{u_2(\lambda_1)}{-2\varepsilon H-\lambda_1}$
into \eqref{4.1-case2-i-j--3.1case1}, and then combining \eqref{4.1-case2--proper-td--3.1} and \eqref{4.1-case2--k--3.1case1} to eliminate
$u_2u_2(H)$ and $u_2u_2(\lambda_1)$, we derive

\begin{equation}\label{4.1-case2--1     3.1case1}
\begin{aligned}
\big(8\Gamma_{1_12}^{1_1}+4\Gamma_{32}^{3}\big)u_2(H)
=&-168\varepsilon\varepsilon_2 H^3+202\varepsilon_2H^2\lambda_1
-64\varepsilon\varepsilon_2H\lambda_1^2+4\varepsilon_2\lambda_1^3\\
&+26c\varepsilon_2 H-2c\varepsilon\varepsilon_2\lambda_1+3\varepsilon_2\lambda H.
\end{aligned}
\end{equation}
Differentiate (\ref{4.1-case2--1     3.1case1}) with $u_2$,
applying \eqref{4.1-case2--proper-td--3.1}, \eqref{4.1-case2-i-j--3.1case1} and \eqref{4.1-case2--k--3.1case1}, we obtain
\begin{equation}\label{4.1-case2--2     3.1case1}
\begin{aligned}
f_1(H,\lambda_1)\Gamma_{1_12}^{1_1}+g_1(H,\lambda_1)\Gamma_{32}^{3}=h_1(H,\lambda_1)u_2(H),
\end{aligned}
\end{equation}
where
\begin{equation*}
\begin{aligned}
&f_1(H,\lambda_1)=1228\varepsilon H ^3 + 136 \varepsilon
H  \lambda_1 ^2 - 852  H ^2 \lambda_1  + 4 c \varepsilon  \lambda_1  - 82 c  H  - 17 \lambda H ,\\
&g_1(H,\lambda_1)=496 \varepsilon H ^3 + 152 \varepsilon  H
\lambda_1 ^2 - 500  H ^2 \lambda_1  - 8  \lambda_1 ^3 + 4 c \varepsilon \lambda_1  - 52 c H  - 10  \lambda H ,\\
&h_1(H,\lambda_1)=600 H ^2 \varepsilon  + 88 \varepsilon \lambda_1 ^2 - 500  H  \lambda_1  - 50 c  - 3 \lambda.
\end{aligned}
\end{equation*}

By differentiating $2\lambda_1+\lambda_3=6\varepsilon H$ along $u_2$, using \eqref{4.1-case2--3--3.1case1}, we know
\begin{equation}\label{4.1-case2--3nu_n(H)  3.1case1}
3\varepsilon u_2(H)=
-(2\varepsilon H+\lambda_1)\Gamma_{1_12}^{1_1}
+(-4\varepsilon H+\lambda_1)\Gamma_{32}^3.
\end{equation}
Putting (\ref{4.1-case2--3nu_n(H)  3.1case1}) into (\ref{4.1-case2--2     3.1case1}), then
\begin{equation}\label{4.1-case2--3     3.1case1}
\begin{aligned}
f_2(H,\lambda_1)\Gamma_{1_12}^{1_1}+g_2(H,\lambda_1)\Gamma_{32}^{3}=0,
\end{aligned}
\end{equation}
where
\begin{equation*}
\begin{aligned}
f_2(H,\lambda_1)=&2484 H ^3 -2156 H ^2 \varepsilon \lambda_1  - 88 \varepsilon \lambda_1 ^3 + 732 H  \lambda_1 ^2 \\&- (146 c + 45 \lambda) H  \varepsilon + (62 c + 3 \lambda)\lambda_1,\\
g_2(H,\lambda_1)=&16(93-450 \varepsilon) H ^3 - 300(5 \varepsilon - 26) \lambda_1  H ^2+ 12 \lambda_1  c \\
 &- 12(213 \varepsilon  -38) H  \lambda_1 ^2- 6(26c \varepsilon+ 5\varepsilon \lambda - 100   c\\& - 6\lambda) H- 24( \varepsilon- 11) \lambda_1 ^3 - 3(50  c  + 3 \lambda)\lambda_1 \varepsilon.
\end{aligned}
\end{equation*}

Multiply both sides of (\ref{4.1-case2--3nu_n(H)  3.1case1}) by $8\Gamma_{1_12}^{1_1}+4\Gamma_{32}^{3}$, applying (\ref{4.1-case2--k--3.1case1}) and (\ref{4.1-case2--1     3.1case1}), we have
\begin{equation}\label{4.1-case2--4     3.1case1}
\begin{aligned}
(16\varepsilon H + 8\lambda_1)(\Gamma_{1_12}^{1_1})^2+(16\varepsilon H - 4\lambda_1)(\Gamma_{32}^{3})^2=h_3(H,\lambda_1),
\end{aligned}
\end{equation}
where
\begin{equation*}
\begin{aligned}
h_3(H,\lambda_1)=&- 366\varepsilon\varepsilon_2H^2\lambda_1 - 4\varepsilon\varepsilon_2\lambda_1^3 + 504\varepsilon_2H^3 +88\varepsilon_2H\lambda_1^2\\&
 - 38\varepsilon c\varepsilon_2H - 9\varepsilon\varepsilon_2\lambda H+ 2c\varepsilon_2\lambda_1.
\end{aligned}
\end{equation*}

By use of (\ref{4.1-case2--k--3.1case1}), \eqref{4.1-case2--3     3.1case1} and \eqref{4.1-case2--4     3.1case1}, we can get a polynomial equation about $H$ and $\lambda_1$
\begin{equation*}
\begin{aligned}
  &(109824 - 1053952 \varepsilon)\lambda_1^9  + (21186688 H -
    2300928 H \varepsilon)\lambda_1^8  + \sum_{i=0}^7r_i\lambda_1^i=0,
\end{aligned}
\end{equation*}
where $r_i$, $0\leq i\leq 7$ are polynomials about $H$. Acting on the above equation by $u_3$ twice, and using
\eqref{4.1-case2-i-j--3.1case1}, (\ref{4.1-case2--k--3.1case1}),
(\ref{4.1-case2--3nu_n(H)  3.1case1}) and (\ref{4.1-case2--3     3.1case1}),
we obtain another algebraic equation of $H$ and $\lambda_1$
\begin{equation*}
\begin{aligned}
(-13895028523008\varepsilon & +
    110378156212224)\lambda_1^{17} - (4577376380387328\varepsilon\\&
    -672576751853568)H\lambda_1^{16}+ \sum_{j=0}^{15}s_j\lambda_1^j=0,
\end{aligned}
\end{equation*}
where $s_j$, $0\leq j\leq 15$ are polynomials about $H$.
Eliminate $\lambda_1$ from the above two equations,
we derive an algebraic equation for $H$ with constant coefficients.
So, $H$ must be a constant, a contradiction.$\hfill\square$


\subsection{The shape operator has the form (I\!I\!I)}

\vskip.2cm
\noindent
{\bf Proposition 3.7}\quad \emph{Let $M^4_r$ be a nondegenerate hypersurface of $N^{5}_s(c)$ with proper mean curvature vector field. Suppose that the shape operator $A$ of $M^4_r$ has the form (I\!I\!I),
then $M^4_r$ has constant mean curvature.}

\vskip.1cm

{\bf Proof}\quad Assume that $H$ is not a constant, then $\nabla H$ is an eigenvector of $A$ with corresponding eigenvalue $-2\varepsilon H$.
Let $\nabla_{u_B}u_C=\Gamma^D_{BC}u_D$, $B, C=1_1, 1_2, 2_1, 2_2$, then
\begin{equation*}
\Gamma_{Di_a}^{j_b}=-\varepsilon_i\varepsilon_j\Gamma_{Dj_{3-b}}^{i_{3-a}},\ i, j, a, b=1, 2,
\end{equation*}
with $D=1_1, 1_2, 2_1, 2_2$. Observe the form (I\!I\!I) of $A$, we know eigenvector $\nabla H$ is in the direction $u_{1_2}$ or $u_{2_2}$. Without loss of generality, we suppose $\nabla H$ is in the direction $u_{1_2}$, then $\lambda_1=-2\varepsilon H$ and $\lambda_2=4\varepsilon H$.

It follows
that
\begin{equation*}
u_{1_{1}}(H)\neq0,\ u_{1_2}(H)=u_{2_1}(H)=u_{2_2}(H)=0,
\end{equation*}
and
\begin{equation*}
\Gamma^{1_{1}}_{BC}=\Gamma^{1_{1}}_{CB},\quad B, C\neq 1_1.
\end{equation*}

We deduce from Codazzi equation that
       \begin{equation*}
       \Gamma_{1_12_2}^{1_1}=\Gamma_{2_22_2}^{1_1}=\Gamma_{2_11_2}^{2_1}=0.
      \end{equation*}
Using Gauss equation for $\langle R(u_{1_{1}}, u_{2_1})u_{1_2},u_{2_2}\rangle$,
combining the above, we get
\begin{equation*}
c-8\varepsilon H^2=0,
\end{equation*}
which tells us $H$ is a constant, a contradiction.
$\hfill\square$

\subsection{The shape operator has the form (I\!V)}

\vskip.2cm
\noindent
{\bf Proposition 3.8}\quad \emph{Let $M^4_r$ be a nondegenerate hypersurface of $N^{5}_s(c)$ with proper mean curvature vector field. Suppose that the shape operator $A$ of $M^4_r$ has the form (I\!V),
then $M^4_r$ has constant mean curvature.}

\vskip.1cm
{\bf Proof}\quad Denote $u_2=u_{2_1}$, and let $\nabla_{u_B}u_C=\Gamma^D_{BC}u_D$, $B, C=1_1, 1_2, 1_3, 2$.
We easily find
\begin{equation}\label{3.2---compatibility-eq1--4}
\Gamma_{D1_1}^{1_3}=\Gamma_{D1_3}^{1_1}=\Gamma_{D1_2}^{1_2}=\Gamma_{D2}^{2}=0,
\end{equation}
and
\begin{equation}\label{3.2---compatibility-eq2--4}
\begin{aligned}
\Gamma_{D1_a}^{1_b}=-\Gamma_{D1_{4-b}}^{1_{4-a}},\ \Gamma_{D1_a}^{2}=-\varepsilon_1\varepsilon_2\Gamma_{D2}^{1_{4-a}},
\end{aligned}
\end{equation}
where $a, b=1, 2, 3$.
Assume that $H$ is not a constant, then $\nabla H$ is an eigenvector of $A$,
is light-like or not.

\vskip.3cm
\emph{Case 1:\ $\nabla H$ is light-like.}

In this case, $\nabla H$ is in the direction $u_{1_3}$ and $\lambda_1=-2\varepsilon H$, $\lambda_2=10\varepsilon H$.
It follows that
\begin{equation*}
u_{1_{1}}(H)\neq0,\ u_{1_2}(H)=u_{1_3}(H)=u_2(H)=0.
\end{equation*}
and
\begin{equation*}
\Gamma^{1_{1}}_{BC}=\Gamma^{1_{1}}_{CB},\quad B, C\neq 1_1.
\end{equation*}

Applying Codazzi equation, we conclude
\begin{equation*}
\Gamma_{1_12}^{1_1}=\Gamma_{1_22}^{1_1}=\Gamma_{21_2}^{2}=\Gamma_{21_3}^{2}=0.
\end{equation*}
Calculate $\langle R(u_{1_{1}}, u_2)u_{1_3},u_2\rangle$ by Gauss equation,
combining the above equations,
we get
\begin{equation*}
c\varepsilon_1-
20\varepsilon\varepsilon_1H^2=0,
\end{equation*}
a contradiction.

\vskip.3cm
\emph{Case 2:\ $\nabla H$ is not light-like.}

We know $\nabla H$ is in the direction $u_2$, $\lambda_2=-2\varepsilon H$ and $\lambda_1=2\varepsilon H$.
So,
\begin{equation}\label{3.2-case2--B(H)--4case2}
u_2(H)\neq0,\ u_B(H)=0, B\neq 2.
\end{equation}
and
\begin{equation}\label{3.2-case2--symmetry--4case2}
\Gamma^2_{BC}=\Gamma^2_{CB},\quad B, C\neq 2.
\end{equation}

We obtain from Codazzi equation that
\begin{equation*}
\Gamma_{21_1}^{2}=\Gamma_{21_2}^{2}=\Gamma_{21_3}^{2}=\Gamma^{2}_{1_{3}1_2}=\Gamma^{2}_{1_{3}1_3}= \Gamma^{1_{1}}_{21_{2}}=\Gamma^{1_{2}}_{21_{3}}=0,
\end{equation*}
       \begin{equation*}
       W:=\Gamma^{1_1}_{1_12}=\Gamma^{1_2}_{1_22}=\Gamma^{1_3}_{1_32},\ \  u_2(H)=-2HW,
       \end{equation*}
Using Gauss equation for $\langle R(u_{2}, u_{1_3})u_{2},u_{1_1}\rangle$, combining the above equations, it gives
\begin{equation}\label{3.2-case2--i--3.1case1}
u_{2}(W)+W^2=4\varepsilon_1 H^2-c\varepsilon_1.
\end{equation}

Let $e_1=\frac{u_{1_1}-u_{1_3}}{\sqrt{2}}$, $e_2=u_{1_2}$, and $e_3=\frac{u_{1_1}+u_{1_3}}{\sqrt{2}}$,
then $\mathfrak{E}=\{e_1, e_2, e_3, u_2\}$ is an orthonormal basis of $T_pM^4_r$, and
(\ref{proper-eq2}) can be rewriten as
\begin{equation*}
\begin{aligned}
2Wu_2(H)
+u_2u_2(H)-16\varepsilon_1\varepsilon H^3+\lambda\varepsilon_1H=0.
\end{aligned}
\end{equation*}
Put $u_2(H)=-2HW$ into the above equation, combining \eqref{3.2-case2--i--3.1case1}, we have
\begin{equation}\label{3.2-case2--proper-td1-3.2.1--bu1}
\begin{aligned}
2HW^2-8(1+2\varepsilon)\varepsilon_1H^3+2c\varepsilon_1H+\lambda\varepsilon_1H=0.
\end{aligned}
\end{equation}
Differentiate \eqref{3.2-case2--proper-td1-3.2.1--bu1} along $u_2$, using $u_2(H)=-2HW$ and \eqref{3.2-case2--i--3.1case1}, we obtain
$$
-2W[4HW^2-16(2+3\varepsilon)\varepsilon_1H^3+4c\varepsilon_1H+\lambda\varepsilon_1H]=0,
$$
which together with \eqref{3.2-case2--proper-td1-3.2.1--bu1} implies
$$
2W[16(1+\varepsilon)\varepsilon_1H^3+\lambda\varepsilon_1H]=0.
$$
Since $u_2(H)\neq 0$, we know $W\neq 0$. The above equation implies $H$ is a constant, a contradiction.
$\hfill\square$

\subsection{The shape operator has the form (V)}

\vskip.2cm
\noindent
{\bf Proposition 3.9}\quad \emph{Let $M^4_r$ be a nondegenerate hypersurface of $N^{5}_s(c)$ with proper mean curvature vector field. Suppose that the shape operator $A$ of $M^4_r$ has the form (V),
then $M^4_r$ has constant mean curvature.}

\vskip.1cm
For the form (V), 
the equations \eqref{5---codazzi-eq--2-c2} and \eqref{5---Gauss-eq--1-c2} deduced from Codazzi equation and Gauss equation are very complicated, compared with the equations for the forms (I), (I\!I), (I\!I\!I) and (I\!V).
However, by constructing creatively the terms $b_k$, $c_k$ (see \eqref{bkck}), and letting $f_k=(\Gamma_{21}^2)^k+2b_k$, with $k=1, \cdots, 5$, we provide an opportunity to complete the proof similarly with the form (I).

\vskip.1cm
{\bf Proof}\quad
Denote $u_1=u_{1_1}$, $u_2=u_{2_1}$, $u_{\bar{3}}=u_{\bar{3}_1}$ and $u_{\tilde{3}}=u_{\tilde{3}_1}$. Let $\nabla_{u_B}u_C=\Gamma_{BC}^Du_D$, $B, C=1, 2, \bar{3}, \tilde{3}$, we have
\begin{equation*}
\Gamma_{BD}^D=0,\ \ \Gamma_{B2}^1=-\varepsilon_1\varepsilon_2\Gamma_{B1}^2,\ \ \Gamma_{Bi}^{\bar{3}}=-\varepsilon_i\Gamma_{B\bar{3}}^i,\ \ \Gamma_{Bi}^{\tilde{3}}=\varepsilon_i\Gamma_{B\tilde{3}}^i,\ \
\Gamma_{B\bar{3}}^{\tilde{3}}=\Gamma_{B\tilde{3}}^{\bar{3}},
\end{equation*}
with $B, D=1, 2, \bar{3}, \tilde{3}$ and $i=1, 2$.

Assume that $H$ is not a constant, then there exists a neighbourhood $U_p$ of $p$  such that $H\neq 0$ and $\nabla H\neq0$.
It follows from \eqref{proper-eq1} that $\nabla H$ is an eigenvector of $A$, and is in the direction $u_1$ or $u_2$.
Without loss of generality, we suppose $\nabla H$ is in
the direction $u_{1}$, then $\lambda_1=-2\varepsilon H$,
\begin{equation*}
u_1(H)\neq 0,\ u_2(H)=u_{\bar{3}}(H)=u_{\tilde{3}}(H)=0.
\end{equation*}
and
\begin{equation*}
\Gamma_{BD}^1=\Gamma_{DB}^1.
\end{equation*}
The equation \eqref{proper-eq2} gives
\begin{equation}\label{5--proper-td-bu-1}
u_1u_1(H)+(\Gamma_{21}^{2}+\Gamma_{\bar{3}1}^{\bar{3}}+\Gamma_{\tilde{3}1}^{\tilde{3}})u_1(H)-\varepsilon \varepsilon_1H\text{tr}A^2+\varepsilon_1\lambda H=0.
\end{equation}

We get from \eqref{Codazzi-td-eq}, with $(X, Y, Z)=(u_1, u_B, u_1), (u_B, u_1, u_B), (u_B, u_D, u_1)$, $B, D=2, \bar{3}, \tilde{3}$ that
\begin{equation}\label{5---codazzi-eq--1-c2}
\Gamma_{1B}^1=\Gamma_{2\bar{3}}^1=\Gamma_{2\tilde{3}}^1=0,\ \ \Gamma_{\bar{3}1}^{\bar{3}}=\Gamma_{\tilde{3}1}^{\tilde{3}},\ \ \text{with}\ B=2, \bar{3}, \tilde{3},
\end{equation}
and
\begin{equation}\label{5---codazzi-eq--2-c2}
\begin{cases}
u_1(\lambda_2)=(\lambda_1-\lambda_2)\Gamma_{21}^2,\\
u_1(\gamma_3)=(\lambda_1-\gamma_3)\Gamma_{\bar{3}1}^{\bar{3}}+\tau_3\Gamma_{\tilde{3}1}^{\bar{3}},\\
u_1(\tau_3)=(\lambda_1-\gamma_3)\Gamma_{\tilde{3}1}^{\bar{3}}-\tau_3\Gamma_{\tilde{3}1}^{\tilde{3}}.
\end{cases}
\end{equation}
If $\lambda_2=\lambda_1$, then the above equations imply $u_1(H)=0$, which contradicts with $u_1(H)\neq0$.
So, $\lambda_2\neq\lambda_1$.
Using Gauss equation for $\langle R(u_1, u_B)u_B, u_1\rangle$ and $\langle R(u_1, u_{\bar{3}})u_{\tilde{3}}, u_1\rangle$, with $B=2, \bar{3}, \tilde{3}$, we have
\begin{equation}\label{5---Gauss-eq--1-c2}
\begin{cases}
u_1(\Gamma_{21}^2)=-(\Gamma_{21}^2)^2+(2H\lambda_2-c)\varepsilon_1,\\
u_1(\Gamma_{\bar{3}1}^{\bar{3}})=-(\Gamma_{\bar{3}1}^{\bar{3}})^2+(\Gamma_{\tilde{3}1}^{\bar{3}})^2+(2H\gamma_3-c)\varepsilon_1,\\
u_1(\Gamma_{\tilde{3}1}^{\bar{3}})=-2\Gamma_{\bar{3}1}^{\bar{3}}\Gamma_{\tilde{3}1}^{\bar{3}}+2H\tau_3\varepsilon_1.
\end{cases}
\end{equation}

{\bf Construct the terms $b_k$ and $c_k$, with $k=1, \cdots, 5$ as follows
\begin{equation}\label{bkck}
\begin{cases}
b_1=\Gamma_{\bar{3}1}^{\bar{3}},\\
c_1=\Gamma_{\tilde{3}1}^{\bar{3}},\\
b_l=b_1b_{l-1}-c_1c_{l-1},\ l=2, 3, 4, 5\\
c_l=b_1c_{l-1}+c_1b_{l-1},\ l=2, 3, 4, 5.
\end{cases}
\end{equation}}
Set $f_k=(\Gamma_{21}^2)^k+2b_k$, with $k=1, \cdots, 5$, we will find \eqref{5---lemma5.2--1} and $u_B(f_1)=0$ also hold, with $B=2, \bar{3}, \tilde{3}$. And then, we can conclude that $u_B(\lambda_2)=u_B(\gamma_3)=u_B(\tau_3)=0$, $B=2, \bar{3}, \tilde{3}$, similarly with Lemma 3.5, which will play an important role in the following process of the proof.

\vskip.2cm
\noindent
{\bf Lemma 3.10}\quad \emph{\eqref{5---lemma5.2--1} holds.}

{\bf Proof}\quad As $\lambda_1=-2\varepsilon H$, the equation \eqref{trA--td---sec2.3} gives $\lambda_2+2\gamma_3=6\varepsilon H$. Differentiate $f_k=(\Gamma_{21}^2)^k+2b_k$ along $u_1$, with $k=1, 2, 3, 4$, using \eqref{5---Gauss-eq--1-c2}, 
we have
\begin{equation}\label{5---lemma5.2-proof-1-c2}
f_{k+1}=-\frac{1}{k}u_1(f_k)-c\varepsilon_1f_{k-1}+2\varepsilon_1HQ_{k-1},
\end{equation}
where $f_0=1$ and $Q_{k-1}=\lambda_2(\Gamma_{21}^2)^{k-1}+2\gamma_3b_{k-1}-2\tau_3c_{k-1}$.
Here, $b_0=1$ and $c_0=0$.

Take the $k$-th derivatives of $\lambda_2+2\gamma_3=6\varepsilon H$ along $u_1$, with $k=1, 2, 3$, combining \eqref{5---codazzi-eq--2-c2} and \eqref{5---Gauss-eq--1-c2}, we obain
\begin{equation}\label{5---lemma5.2-proof-2-c2}
Q_1=-6\varepsilon u_1(H)-2\varepsilon Hf_1,
\end{equation}
\begin{equation*}
Q_2=3\varepsilon u_1^{(2)}(H)+\varepsilon Hu_1(f_1)+\varepsilon f_1u_1(H)-\varepsilon Hf_2-4\varepsilon_1H^3-3\varepsilon\varepsilon_1Hc+\varepsilon_1H\text{tr}A^2,
\end{equation*}
and
\begin{equation*}
\begin{aligned}
3Q_3=&\varepsilon Hu_1(f_2)-3\varepsilon u_1^{(3)}(H)+\varepsilon f_2u_1(H)-2\varepsilon u_1(H)u_1(f_1)\\
&+12\varepsilon_1H^2u_1(H)-\varepsilon_1u_1(H)\text{tr}A^2-\varepsilon f_1u_1^{(2)}(H)+15\varepsilon\varepsilon_1cu_1(H)\\
&-\varepsilon Hu_1^{(2)}(f_1)+4\varepsilon\varepsilon_1 cHf_1-2\varepsilon Hf_3-\varepsilon_1Hu_1(\text{tr}A^2)\\&+4\varepsilon_1H(\lambda_2^2\Gamma_{21}^2+2(\gamma_3^2-\tau_3^2)\Gamma_{\bar{3}1}^{\bar{3}}-4\gamma_3\tau_3\Gamma_{\tilde{3}1}^{\bar{3}}).
\end{aligned}
\end{equation*}

Multiply the first, second and third equation in \eqref{5---codazzi-eq--2-c2} by $\lambda_2$, $\gamma_3$ and $\tau_3$, respectively, and then we get from these results that
\begin{equation*}
\lambda_2^2\Gamma_{21}^2+2(\gamma_3^2-\tau_3^2)\Gamma_{\bar{3}1}^{\bar{3}}-4\gamma_3\tau_3\Gamma_{\tilde{3}1}^{\bar{3}}=-2\varepsilon HQ_1-\frac{1}{2}u_1(\text{tr}A^2)+4Hu_1(H).
\end{equation*}
Applying \eqref{5--proper-td-bu-1} and the above equations,
we deduce from \eqref{5---lemma5.2-proof-1-c2} that \eqref{5---lemma5.2--1} holds.
$\hfill\square$

By calculation, we easily find \eqref{5---fkrelations} also holds. Then, follow the process of the proof for Lemma 3.4 (replacing the index $i\in\{2, 3, 4\}$ with $B\in\{2, \bar{3}, \tilde{3}\}$), we conclude that $u_B(f_1)=0$, with $B=2, \bar{3}, \tilde{3}$.

\vskip.2cm
\noindent
{\bf Lemma 3.11}\quad \emph{We have $$
u_B(\Gamma_{D1}^D)=u_B(\Gamma_{\tilde{3}1}^{\bar{3}})=u_B(\lambda_2)=u_B(\gamma_3)=u_B(\tau_3)=0,
$$
with $B, D=2, \bar{3}, \tilde{3}$.}

{\bf Proof}\quad Since $u_B(f_1)=u_B(H)=0$, $B=2, \bar{3}, \tilde{3}$, we conclude
$$
u_Bu_1^{(k)}(f_1)=u_Bu_1^{(k)}(H)=0,\ B=2, \bar{3}, \tilde{3},\ k=1, 2, 3.
$$
Observe the expressions of $f_k$, $k=2, 3$ in \eqref{5---lemma5.2--1}, combining the above equation,
we have
\begin{equation*}
\begin{aligned}
u_B(f_k)=0,\  B=2, \bar{3}, \tilde{3},\ k=1, 2, 3,
\end{aligned}
\end{equation*}
which gives
\begin{equation*}
\begin{aligned}
(\Gamma_{21}^{2})^{k-1}u_B(\Gamma_{21}^{2})+2b_{k-1}u_B(\Gamma_{\bar{3}1}^{\bar{3}})-2c_{k-1}u_B(\Gamma_{\tilde{3}1}^{\bar{3}})=0, k=1, 2, 3.
\end{aligned}
\end{equation*}
Here, $b_0=1$ and $c_0=0$.
The coefficient determinant of this system
\begin{equation*}
\begin{aligned}
\left|
  \begin{array}{ccc}
     1   &2      &0\\
      \Gamma_{21}^{2}   & 2\Gamma_{\bar{3}1}^{\bar{3}} & -2\Gamma_{\tilde{3}1}^{\bar{3}}\\
     (\Gamma_{21}^{2})^2&     2((\Gamma_{\bar{3}1}^{\bar{3}})^2-(\Gamma_{\tilde{3}1}^{\bar{3}})^2)& -4\Gamma_{\bar{3}1}^{\bar{3}}\Gamma_{\tilde{3}1}^{\bar{3}}
  \end{array}
\right|=-4\Gamma_{\tilde{3}1}^{\bar{3}}[(\Gamma_{21}^2-\Gamma_{\bar{3}1}^{\bar{3}})^2+(\Gamma_{\tilde{3}1}^{\bar{3}})^2].
\end{aligned}
\end{equation*}
As $\tau_3\neq 0$, the third equation of \eqref{5---Gauss-eq--1-c2} implies that $\Gamma_{\tilde{3}1}^{\bar{3}}\neq 0$. So, the above coefficient determinant is not equal to zero, and then
\begin{equation*}
u_B(\Gamma_{21}^{2})=u_B(\Gamma_{\bar{3}1}^{\bar{3}})=u_B(\Gamma_{\tilde{3}1}^{\bar{3}})=0,\ B=2, \bar{3}, \tilde{3}.
\end{equation*}
Differentiate the equations in \eqref{5---Gauss-eq--1-c2} along $u_B$, $B=2, \bar{3}, \tilde{3}$, combining the above equation,
we get
$$
u_B(\lambda_2)=u_B(\gamma_3)=u_B(\tau_3)=0,\ B=2, \bar{3}, \tilde{3}.
$$
$\hfill\square$

Now, we continue the proof of Proposition 3.9.

Applying Lemma 3.11, we obtain from the equation \eqref{Codazzi-td-eq} for $(X, Y, Z)=(u_B, u_D, u_B), (u_2, u_{\bar{3}}, u_{\tilde{3}}), (u_2, u_{\tilde{3}}, u_{\bar{3}})$, with $B, D=2, \bar{3}, \tilde{3}$ that
\begin{equation}\label{codazzi}
\begin{cases}
(\gamma_3-\lambda_2)\Gamma_{2\bar{3}}^2-\tau_3\Gamma_{2\tilde{3}}^2=0,\\
(\gamma_3-\lambda_2)\Gamma_{2\tilde{3}}^2+\tau_3\Gamma_{2\bar{3}}^2=0,\\
-2\tau_3\Gamma_{2\tilde{3}}^{\bar{3}}=(\lambda_2-\gamma_3)\Gamma_{\bar{3}2}^{\bar{3}}-\tau_3\Gamma_{\bar{3}2}^{\tilde{3}},\\
2\tau_3\Gamma_{2\tilde{3}}^{\bar{3}}=(\lambda_2-\gamma_3)\Gamma_{\tilde{3}2}^{\tilde{3}}+\tau_3\Gamma_{\tilde{3}2}^{\bar{3}},\\
(\lambda_2-\gamma_3)\Gamma_{\bar{3}2}^{\tilde{3}}+\tau_3\Gamma_{\bar{3}2}^{\bar{3}}=0,\\
(\lambda_2-\gamma_3)\Gamma_{\tilde{3}2}^{\bar{3}}-\tau_3\Gamma_{\tilde{3}2}^{\tilde{3}}=0,\\
\Gamma_{\tilde{3}\bar{3}}^{\tilde{3}}=\Gamma_{\bar{3}\tilde{3}}^{\bar{3}}=0,
\end{cases}
\end{equation}
which gives that
\begin{equation*}
\Gamma_{2\bar{3}}^2=\Gamma_{2\tilde{3}}^2=\Gamma_{\bar{3}\tilde{3}}^{\bar{3}}=\Gamma_{\tilde{3}\bar{3}}^{\tilde{3}}=\Gamma_{\bar{3}2}^{\bar{3}}+\Gamma_{\tilde{3}2}^{\tilde{3}}=\Gamma_{\tilde{3}2}^{\bar{3}}-\Gamma_{\bar{3}2}^{\tilde{3}}=0,
\end{equation*}
and
\begin{equation}\label{5---codazzi-td--2-c2}
(\Gamma_{\bar{3}2}^{\bar{3}})^2+(\Gamma_{\bar{3}2}^{\tilde{3}})^2-2\Gamma_{\bar{3}2}^{\tilde{3}}\Gamma_{2\bar{3}}^{\tilde{3}}=0.
\end{equation}

Calculate $\langle R(u_2, u_B)u_2, u_B\rangle$, $\langle R(u_2, u_B)u_2, u_D\rangle$, 
$\langle R(u_B, u_D)u_B, u_D\rangle$ and $\langle R(u_2, u_B)u_D, u_1\rangle$ by Gauss equation, with $B, D=\bar{3}, \tilde{3}$ and $B\neq D$, combining \eqref{5---codazzi-eq--1-c2} and the above two equations, we have
\begin{equation}\label{4zql-x-1--c2}
\Gamma_{\bar{3}2}^{\tilde{3}}(\Gamma_{21}^2-\Gamma_{\bar{3}1}^{\bar{3}})+\Gamma_{\bar{3}2}^{\bar{3}}\Gamma_{\tilde{3}1}^{\bar{3}}=0,
\end{equation}
and
\begin{equation}\label{5---Gauss-eq--3qian--c2}
\begin{cases}
-\varepsilon_1\Gamma_{21}^2\Gamma_{\tilde{3}1}^{\bar{3}}+2\varepsilon_2\Gamma_{\bar{3}2}^{\bar{3}}\Gamma_{2\bar{3}}^{\tilde{3}}=\varepsilon\lambda_2\tau_3,\\
\varepsilon_1\Gamma_{21}^2\Gamma_{\bar{3}1}^{\bar{3}}+2\varepsilon_2\Gamma_{\bar{3}2}^{\tilde{3}}\Gamma_{2\bar{3}}^{\tilde{3}}=-c-\varepsilon\lambda_2\gamma_3,\\
\varepsilon_1[(\Gamma_{\bar{3}1}^{\bar{3}})^2+(\Gamma_{\tilde{3}1}^{\bar{3}})^2]-4\varepsilon_2\Gamma_{\bar{3}2}^{\tilde{3}}\Gamma_{2\bar{3}}^{\tilde{3}}=-c-\varepsilon(\gamma_3^2+\tau_3^2).
\end{cases}
\end{equation}
The second and third equations of \eqref{5---Gauss-eq--3qian--c2} implies
\begin{equation}\label{5---Gauss-eq--3--c2}
\varepsilon_1[(\Gamma_{\bar{3}1}^{\bar{3}})^2+(\Gamma_{\tilde{3}1}^{\bar{3}})^2+2\Gamma_{21}^2\Gamma_{\bar{3}1}^{\bar{3}}]=-3c-2\varepsilon\lambda_2\gamma_3-\varepsilon(\gamma_3^2+\tau_3^2).
\end{equation}

\vskip.2cm
\emph{Case 1:\ $\Gamma_{\bar{3}2}^{\tilde{3}}\neq 0$ at some point in $U_p$.}

\vskip.1cm

We suppose $\Gamma_{\bar{3}2}^{\tilde{3}}\neq 0$ on a neighbourhood  $U_q\subset U_p$, and work on $U_q$. From the fifth equation in \eqref{codazzi} and \eqref{4zql-x-1--c2}, we get
$$
\frac{\Gamma_{21}^2-\Gamma_{\bar{3}1}^{\bar{3}}}{\Gamma_{\tilde{3}1}^{\bar{3}}}=\frac{\lambda_2-\gamma_3}{\tau_3}.
$$
Considering Lemma 3.11, we conclude from the above equation that there exists two smooth functions $\mu$ and $\nu$, with $u_B(\mu)=u_B(\nu)=0$, $B=2, \bar{3}, \tilde{3}$, such that
\begin{equation}\label{5---gamma--4--c2}
\Gamma_{21}^2=\mu\lambda_2+\nu,\ \Gamma_{\bar{3}1}^{\bar{3}}=\mu\gamma_3+\nu,\ \Gamma_{\tilde{3}1}^{\bar{3}}=\mu\tau_3.
\end{equation}
The equation \eqref{5---f1--2} also holds. By use of \eqref{5---codazzi-eq--2-c2} and \eqref{5---Gauss-eq--1-c2}, differentiating the equations in \eqref{5---gamma--4--c2} along $u_1$,
it gives \eqref{5---u1mu-u1nu-1}.
Substitute \eqref{5---gamma--4--c2} into \eqref{5---Gauss-eq--3--c2}, we obtain
\begin{equation}\label{subsection5--bu1--case1}
(\varepsilon_1\mu^2+\varepsilon)(\gamma_3^2+2\lambda_2\gamma_3+\tau_3^2)+2\varepsilon_1\mu\nu(\lambda_2+2\gamma_3)+3\varepsilon_1\nu^2
+3c=0.
\end{equation}
Since
$$
\begin{aligned}
\gamma_3^2+2\lambda_2\gamma_3+\tau_3^2&=\frac{1}{2}[(\lambda_2+2\gamma_3)^2-(\lambda_2^2+2\gamma_3^2-2\tau_3^2)]\\
&=20H^2-\frac{1}{2}\text{tr}A^2,
\end{aligned}
$$
the equation \eqref{subsection5--bu1--case1} can be rewritten as \eqref{5---Gauss-eq--4}.
Put \eqref{5---gamma--4--c2} into \eqref{5---lemma5.2-proof-2-c2}, combining \eqref{5---f1--2} and \eqref{5---Gauss-eq--4}, we have \eqref{5---u1H--1}.
Follow the process for the case $\Gamma_{23}^{4}\neq 0$ at some point in section 3.1, just correcting \eqref{5--proper-td-1} with \eqref{5--proper-td-bu-1}, one derive a contradiction.

\vskip.2cm
\emph{Case 2:\ $\Gamma_{\bar{3}2}^{\tilde{3}}\equiv0$ on $U_p$.}

\vskip.1cm

It follows from \eqref{5---codazzi-td--2-c2} that $\Gamma_{\bar{3}2}^{\bar{3}}=0$. And then, \eqref{5---Gauss-eq--3qian--c2} can be simplified into
\begin{equation}\label{5-II--Gauss-eq--3qian--c2}
\begin{cases}
\varepsilon_1\Gamma_{21}^2\Gamma_{\tilde{3}1}^{\bar{3}}=-\varepsilon\lambda_2\tau_3,\\
\varepsilon_1\Gamma_{21}^2\Gamma_{\bar{3}1}^{\bar{3}}=-c-\varepsilon\lambda_2\gamma_3,\\
\varepsilon_1[(\Gamma_{\bar{3}1}^{\bar{3}})^2+(\Gamma_{\tilde{3}1}^{\bar{3}})^2]=-c-\varepsilon(\gamma_3^2+\tau_3^2).
\end{cases}
\end{equation}

\vskip.2cm
\emph{Subcase (1):\ $\Gamma_{21}^{2}\equiv0$ on $U_p$.}

\vskip.1cm

The first and second equations of \eqref{5-II--Gauss-eq--3qian--c2} tells us that $\lambda_2=0$ and $c+\varepsilon\lambda_2\gamma_3=0$.
If $c\neq 0$, then $c+\varepsilon\lambda_2\gamma_3\neq 0$, 
a contradiction.
In the following, we suppose $c=0$.
Differentiate both sides of $\gamma_3=3\varepsilon H$ along $u_1$, using \eqref{5---codazzi-eq--2-c2}, we get
\begin{equation}\label{5---II.2--u1H--c2}
3\varepsilon u_1(H)=(-2\varepsilon H-\gamma_3)\Gamma_{\bar{3}1}^{\bar{3}}+\tau_3\Gamma_{\tilde{3}1}^{\bar{3}}.
\end{equation}
Applying \eqref{5-II--Gauss-eq--3qian--c2} and \eqref{5---II.2--u1H--c2}, we can rewrite \eqref{5--proper-td-bu-1} as
\begin{equation}\label{5---II.2--propertdx2--c2}
\begin{aligned}
4[-5\varepsilon H((\Gamma_{\bar{3}1}^{\bar{3}})^2-(\Gamma_{\tilde{3}1}^{\bar{3}})^2)+2\tau_3\Gamma_{\bar{3}1}^{\bar{3}}\Gamma_{\tilde{3}1}^{\bar{3}}]-34\varepsilon_1H\tau_3^2
+162\varepsilon_1H^3-9\varepsilon\varepsilon_1\lambda H=0.
\end{aligned}
\end{equation}
By differentiating \eqref{5---II.2--propertdx2--c2} along $u_1$, using \eqref{5---codazzi-eq--2-c2} and \eqref{5---Gauss-eq--1-c2},
and combining \eqref{5-II--Gauss-eq--3qian--c2} and \eqref{5---II.2--propertdx2--c2}, we derive
\begin{equation}\label{5---II.2--upropertdx2--c2}
\begin{aligned}
K_1\Gamma_{\bar{3}1}^{\bar{3}} + K_2\Gamma_{\tilde{3}1}^{\bar{3}}=0,
\end{aligned}
\end{equation}
where
$$
\begin{cases}
K_1=1710H^3 - 78H\tau_3^2 - 135\varepsilon\lambda H,\\
K_2=1494H^2\tau_3\varepsilon - 78\tau_3^3\varepsilon - 9\tau_3\lambda.
\end{cases}
$$
As \eqref{5-II--Gauss-eq--3qian--c2}, we get from \eqref{5---II.2--upropertdx2--c2} that
\begin{equation*}
\begin{cases}
(K_1^2+K_2^2)\Gamma_{\bar{3}1}^{\bar{3}}\Gamma_{\tilde{3}1}^{\bar{3}}-\varepsilon\varepsilon_1(9H^2+\tau_3^2)K_1K_2=0,\\
K_1K_2[(\Gamma_{\bar{3}1}^{\bar{3}})^2-(\Gamma_{\tilde{3}1}^{\bar{3}})^2]+(K_2^2-K_1^2)\Gamma_{\bar{3}1}^{\bar{3}}\Gamma_{\tilde{3}1}^{\bar{3}}=0.
\end{cases}
\end{equation*}
Eliminate $\Gamma_{\bar{3}1}^{\bar{3}}$ and $\Gamma_{\tilde{3}1}^{\bar{3}}$ from \eqref{5---II.2--propertdx2--c2} and the above equations, one derives
\begin{equation}\label{5---II.2--lambda23H---1--c2}
\begin{aligned}
&2741856\varepsilon \tau_3^{12} + c_1\tau_3 ^{10} +c_2\tau_3
^8 + c_3 \tau_3 ^6 + c_4\tau_3 ^4 +c_5\tau_3 ^2+c_6=0,
\end{aligned}
\end{equation}
where $c_i$, $1\leq i\leq 6$, are polynomials about $H$.
Acting on \eqref{5---II.2--lambda23H---1--c2} by $u_1$, we have
\begin{equation*}
\begin{aligned}
L_1\Gamma_{\bar{3}1}^{\bar{3}}+L_2\Gamma_{\tilde{3}1}^{\bar{3}}=0,
\end{aligned}
\end{equation*}
where $L_1$ and $L_2$ are polynomials about $H$ and $\tau_3$.
It follows from \eqref{5---II.2--upropertdx2--c2} and the above equation that
$$
K_1L_2-K_2L_1=0.
$$
Eliminating $\tau_3$ from \eqref{5---II.2--lambda23H---1--c2} and the above equation, we get a polynomial equation of degree $92$ for $H$, which implies $H$ is a constant, a contradiction.

\vskip.2cm
\emph{Subcase (2):\ $\Gamma_{21}^{2}\neq 0$ at some point in $U_p$.}

\vskip.1cm

We suppose $\Gamma_{21}^2\neq 0$ at $q\in U_p$.
If $c+\varepsilon(\gamma_3^2+\tau_3^2)=0$, then \eqref{5-II--Gauss-eq--3qian--c2} implies $\Gamma_{\bar{3}1}^{\bar{3}}=\Gamma_{\tilde{3}1}^{\bar{3}}=0$ and $\lambda_2=c=\gamma_3=\tau_3=0$, a contradiction.
So, $c+\varepsilon(\gamma_3^2+\tau_3^2)\neq 0$. We derive from the equations in \eqref{5-II--Gauss-eq--3qian--c2} that
\begin{equation}\label{5-II--gamma212--c2}
\begin{aligned}
(\Gamma_{21}^{2})^2=-\varepsilon_1\frac{\lambda_2^2(\gamma_3^2+\tau_3^2)+c^2+2c\varepsilon\lambda_2\gamma_3}{c+\varepsilon(\gamma_3^2+\tau_3^2)}.
\end{aligned}
\end{equation}

Let
$$\mu=2\lambda_2\gamma_3+\gamma_3^2+\tau_3^2,\ \
\nu=\lambda_2(\gamma_3^2+\tau_3^2),
$$
then take derivatives of $\mu$, $\nu$ and $H$ along $u_1$, by use of \eqref{5---codazzi-eq--2-c2}, \eqref{5-II--Gauss-eq--3qian--c2} and \eqref{5-II--gamma212--c2}, we get
\begin{equation}\label{5-II--u1H-u1mu-u1nu--c2}
\begin{cases}
6\varepsilon u_1(H)=K(12\varepsilon Hc^2+10H\mu c-3\varepsilon\nu c+48H^2\nu-2\mu\nu),\\
u_1(\mu)=K[(24H^2+2\mu)c^2+(2\mu^2+12\varepsilon H^2\mu-6H\nu)c-3\nu^2+10\varepsilon H\mu\nu],\\
u_1(\nu)=K[(2\varepsilon H\mu+3\nu)c^2+(24\varepsilon H^2\nu+2\varepsilon\mu\nu)c+12\varepsilon H\nu^2],
\end{cases}
\end{equation}
where $K=\frac{\varepsilon_1}{\Gamma_{21}^2[c+\varepsilon(\gamma_3^2+\tau_3^2)]}$.
Applying \eqref{5-II--Gauss-eq--3qian--c2}, \eqref{5-II--gamma212--c2} and \eqref{5-II--u1H-u1mu-u1nu--c2},
we can rewrite \eqref{5--proper-td-bu-1} as \eqref{5--proper-td-2}.
Follow the process for the subcase $\Gamma_{21}^2, \Gamma_{\bar{3}1}^{\bar{3}}, \Gamma_{\tilde{3}1}^{\tilde{3}}\neq 0$ of the case $\Gamma_{2\bar{3}}^{\tilde{3}}=0$ in section 3.1,
and replace the system \eqref{5-II--u1H-u1mu-u1nu} with \eqref{5-II--u1H-u1mu-u1nu--c2}, the equation \eqref{bueq2---sec5} with
\begin{equation*}
\begin{aligned}
-\varepsilon_1(c+(\gamma_3^2+\tau_3^2))^2(\Gamma_{21}^2)^2&=(c+(\gamma_3^2+\tau_3^2))(c^2+\lambda_2^2(\gamma_3^2+\tau_3^2)+2c\lambda_2\gamma_3)\\
&=\nu^2+6cH\nu+c^2\mu+c^3=0,
\end{aligned}
\end{equation*}
we deduce a contradiction.
$\hfill\square$

\vskip.2cm
\noindent
{\bf Remark 3.12}\quad For form (V) of the shape operator, by constructing the terms $b_k$ and $c_k$, $k=1, 2, \cdots$, we find there are many equations similarly with the equations for the form (I) hold, by only changing some terms for the form (I) into the terms about $b_k$ and $c_k$ for the form (V) as follows.
$$
\begin{array}{|c|c|}
  \hline\textbf{The terms for form (I)} &  \textbf{The terms for form (V)}\\
  \hline (\Gamma_{31}^3)^k+(\Gamma_{41}^4)^k & 2b_k\\
 \lambda_3(\Gamma_{31}^3)^k+\lambda_4(\Gamma_{41}^4)^k & 2\gamma_3b_k-2\tau_3c_k\\
  \lambda_3^2(\Gamma_{31}^3)^k+\lambda_4^2(\Gamma_{41}^4)^k & 2(\gamma_3^2-\tau_3^2)b_{k}-4\gamma_3\tau_3c_{k-1} \\
  \hline
\end{array}
$$
We believe that the construction of $b_k$ and $c_k$, $k=1, 2, \cdots$, will provide a new insight for us to study hypersurfaces with imaginary principal curvatures.

\subsection{The shape operator has the form (V\!I)}

\vskip.2cm
\noindent
{\bf Proposition 3.13}\quad \emph{Let $M^4_r$ be a nondegenerate hypersurface of $N^{5}_s(c)$ with proper mean curvature vector field. Suppose that the shape operator $A$ of $M^4_r$ has the form (V\!I),
then $M^4_r$ has constant mean curvature.}

{\bf Proof}\quad Assume that $H$ is not a constant, then the eigenvector $\nabla H$ of $A$ is in the direction $u_{1_2}$ and $\lambda_1=-2\varepsilon H$. As $2\lambda_1+2\gamma_2=4\varepsilon$, we know $\gamma_2=4\varepsilon H$.
Let $u_{\bar{2}}=u_{\bar{2}_1}$, $u_{\tilde{2}}=u_{\tilde{2}_1}$ and $\nabla_{u_B}u_C=\Gamma^D_{BC}u_D$, $B, C=1_1, 1_2, \bar{2}, \tilde{2}$,
then
\begin{equation*}
\Gamma_{D1_1}^{1_2}=\Gamma_{D\bar{2}}^{\bar{2}}=\Gamma_{D\tilde{2}}^{\tilde{2}}=0,
\end{equation*}
and
\begin{equation*}
\begin{aligned}
\Gamma_{D1_1}^{1_1}=-\Gamma_{D1_2}^{1_2},\ \Gamma_{D\bar{2}}^{\tilde{2}}=\Gamma_{D\tilde{2}}^{\bar{2}},\ \Gamma_{D1_a}^{\bar{2}}=-\varepsilon_1\Gamma_{D\bar{2}}^{1_{3-a}},\ \Gamma_{D1_a}^{\tilde{2}}=\varepsilon_1\Gamma_{D\tilde{2}}^{1_{3-a}},\ \ a=1, 2.
\end{aligned}
\end{equation*}
It follows that
\begin{equation*}
u_{1_{1}}(H)\neq0,\ u_{1_2}(H)=u_{\bar{2}}(H)=u_{\tilde{2}}(H)=0.
\end{equation*}
and
\begin{equation*}
\Gamma^{1_{1}}_{BC}=\Gamma^{1_{1}}_{CB},\quad B, C\neq 1_1.
\end{equation*}
From Codazzi equation \eqref{Codazzi-td-eq},
with $(X, Y, Z)=(u_{1_1}, u_{B}, u_{1_2}), (u_{\bar{2}}, u_{\tilde{2}}, u_{1_2}),\\ (u_{1_2}, u_{B}, u_{B}), (u_{1_2}, u_{\bar{2}}, u_{\tilde{2}})$, $B=\bar{2}, \tilde{2}$, we deduce
\begin{equation*}
\Gamma_{\bar{2}\bar{2}}^{1_1}=-\Gamma_{\tilde{2}\tilde{2}}^{1_1},\ \ \Gamma_{1_1\bar{2}}^{1_1}=\Gamma_{1_1\tilde{2}}^{1_1}=\Gamma_{1_2 \bar{2}}^{\tilde{2}}=\Gamma_{1_2 \tilde{2}}^{\bar{2}}=0,
\end{equation*}
and
\begin{equation}\label{4.2---codazzitds3}
\begin{cases}
-6\varepsilon H \Gamma_{\bar{2}1_2}^{\bar{2}}+\tau_2\Gamma_{\tilde{2}1_2}^{\bar{2}}=0,\\
u_{1_2}(\tau_1)=-6\varepsilon H \Gamma_{\tilde{2}1_2}^{\bar{2}}-\tau_2 \Gamma_{\bar{2}1_2}^{\bar{2}}.
\end{cases}
\end{equation}
Using Gauss equation for
$\langle R(u_{1_2}, u_{\bar{2}})u_{\bar{2}}, u_{1_2}\rangle$ and
$\langle R(u_{1_2}, u_{\bar{2}})u_{\tilde{2}}, u_{1_2}\rangle$,
combining the above equations, it gives
\begin{equation}\label{4.2---Gausstds1}
\begin{cases}
u_{1_2}(\Gamma_{\bar{2}1_2}^{\bar{2}})
=\Gamma^{1_2}_{1_21_2}
\Gamma_{\bar{2}1_2}^{\bar{2}}
+(\Gamma_{\tilde{2} 1_2}^{\bar{2}})^2
-(\Gamma_{\bar{2} 1_2}^{\bar{2}})^2,\\
u_{1_2}(\Gamma_{\tilde{2}1_2}^{\bar{2}})
=\Gamma^{1_2}_{1_21_2}
\Gamma_{\tilde{2}1_2}^{\bar{2}}
-2\Gamma_{\bar{2} 1_2}^{\bar{2}}\Gamma_{\tilde{2} 1_2}^{\bar{2}}.
\end{cases}
\end{equation}
Applying \eqref{4.2---codazzitds3} and \eqref{4.2---Gausstds1}, differentiate the first equation in \eqref{4.2---codazzitds3} along $u_{1_2}$,
we obtain
\begin{equation*}
-6\varepsilon H (\Gamma_{\tilde{2}1_2}^{\bar{2}})^2-\tau_2\Gamma_{\bar{2}1_2}^{\bar{2}}\Gamma_{\tilde{2}1_2}^{\bar{2}}=0,
\end{equation*}
which together with the first equation of \eqref{4.2---codazzitds3} implies that
$
\Gamma_{\tilde{2}1_2}^{\bar{2}}=0
$
and
$
H\Gamma_{\bar{2}1_2}^{\bar{2}}=0
$.
Compute $\langle R(u_{1_1}, u_{\bar{2}})u_{1_2}, u_{\tilde{2}}\rangle$
by Gauss equation, we have
$$
\Gamma_{\bar{2}1_2}^{\bar{2}}\Gamma_{1_1\bar{2}}^{\tilde{2}}=2H\tau_2.
$$
Multiply both sides of the above equation by $H$, combining $
H\Gamma_{\bar{2}1_2}^{\bar{2}}=0
$, we conclude that
$
H=0
$,
a contradiction.
$\hfill\square$

\section{Estimates for the constant mean curvature}

In this section, as an application of the conclusion that $M^4_r$ has constant mean curvature,
we will estimate the value of the mean curvature $H$. When $M^4_r$ has one or two distinct principal curvatures, we will compute the value of $H$. It's a pity that we can not get the value of $H$ when
$M^4_r$ has three or four distinct principal curvatures. But for the case that the principal curvatures are all real, we give a value range of the constant mean curvature $H$.

\subsection{When $M^4_r$ has one simple principal curvature}

\vskip.2cm
\noindent
{\bf Theorem 4.1}\quad \emph{Let $M^4_r$ be a nondegenerate hypersurface of $N^{5}_s(c)$ with proper mean curvature vector field, satisfying $\Delta\vec{H}=\lambda\vec{H}$, $\vec{\xi}$ be a unit normal vector field to $M^4_r$, with $\varepsilon=\langle\vec\xi,\vec\xi\rangle=\pm1$. Suppose that $M^4_r$ has one simple principal curvature, then
\begin{itemize}
  \item when $\varepsilon \lambda\leq 0$, we have $H=0$;
  \item when $\varepsilon \lambda>0$, we have $H=0$ or
$
H^2=\frac{\varepsilon\lambda}{4}
$.
\end{itemize}
}

{\bf Proof}\quad Recall subsection 2.3, the imaginary principal curvatures of $M^4_r$ appear in conjugate pairs.
So, under the assumption that $M^4_r$ has one simple principal curvature $\mu$, $\mu$ must be real. Then, we know $t=m$ (the signs $t, m$ see subsection 2.3) and $\lambda_1=\cdots=\lambda_m=\mu$. Thus, equations \eqref{trA--td---sec2.3} can be simplified as $\text{tr}A^2=4\mu^2$ and 
$\mu=\varepsilon H$.
According to Theorem 3.1, the mean curvature $H$ of $M^4_r$ is a constant.
Using this result and $\text{tr}A^2=4\mu^2$, we get from \eqref{proper-eq2} that $H=0$ or $4\mu^2=\varepsilon \lambda$.

When $\varepsilon\lambda\leq 0$, we assume $H\neq 0$. Then $4\mu^2=\varepsilon \lambda\leq 0$, a contradiction.

When $\varepsilon\lambda>0$, suppose that $H\neq 0$, we also have $4\mu^2=\varepsilon \lambda$, which together with $\mu=\varepsilon H$ implies that
$$
H^2=\frac{\varepsilon\lambda}{4}.
$$
$\hfill\square$

\vskip.2cm
\noindent
{\bf Remark 4.2}\quad The result of Theorem 4.1 for $c=0$ has been gotten in \cite{Liu 2017}.

\subsection{When $M^4_r$ has two distinct principal curvatures}

\vskip.2cm
\noindent
{\bf Lemma 4.3}\quad \emph{Let $M^4_r$ be a non-minimal hypersurface of $N^{5}_s(c)$ with proper mean curvature vector field, $\vec{\xi}$ be a unit normal vector field to $M^4_r$, with $\varepsilon=\langle\vec\xi,\vec\xi\rangle=\pm1$. Suppose that $M^4_r$ has two distinct principal curvatures $\mu$ and $\nu$, then
   $$
       c+\varepsilon\mu\nu=0.
       $$}

{\bf Proof}\quad
Since the number of imaginary principal curvatures is even, these two distinct principal curvatures $\mu$ and $\nu$ are all real or all imaginary.

\vskip.2cm
\emph{Case 1:\ $\mu$ and $\nu$ are all real.}

\vskip.1cm
Let $l$ is the multiply of $\mu$, then 
considering the mean curvature $H$ is a constant and $H\neq 0$, we get from \eqref{proper-eq2} and
\eqref{trA--td---sec2.3} that
\begin{equation}\label{trA-trA2-td-tworeal-bu1-sec2.3}
\begin{cases}
l\mu+(4-l)\nu=4\varepsilon H,\\
l\mu^2+(4-l)\nu^2=\varepsilon \lambda,
\end{cases}
\end{equation}
which implies $\mu$ and $\nu$ are all constants.
In the following, by use of that the principal curvatures are all constant, we
discuss separately the possible forms (I), (I\!I), (I\!I\!I) and (I\!V) of $A$ for this case, and derive the result $c+\varepsilon\mu\nu=0$.

\vskip.2cm
\emph{Subcase (1):\ the shape operator has the form (I).}

\vskip.1cm
Without loss of generality, we suppose $\lambda_1=\mu$ and $\lambda_4=\nu$. Let $u_i=u_{i_1}$ and $\nabla_{u_i}u_j=\Gamma_{ij}^ku_k, i, j=1, 2, 3, 4$, then recall subsection 3.1, we know \eqref{5---compatibility-eq1--4} and \eqref{5---codazzi-eq--2} also hold. Since the principal curvatures $\lambda_i$, with $i=1, 2, 3, 4$, are constants, the equations \eqref{5---compatibility-eq1--4} and \eqref{5---codazzi-eq--2} imply that $\Gamma_{ij}^i=0$ if $\lambda_i\neq \lambda_j$, $\Gamma_{i4}^1=\Gamma_{i1}^4=\Gamma_{1i}^4=0$ if $\lambda_i=\mu$, and $\Gamma_{41}^i=\Gamma_{i1}^4=0$ if $\lambda_i=\nu$.
Using Gauss equation for $\langle R(u_1, u_4)u_1, u_4\rangle$, we have
$$
c+\varepsilon\mu\nu=0.
$$

\emph{Subcase (2):\ the shape operator has the form (I\!I).}

\vskip.1cm
We can suppose $\lambda_1=\mu$ and $\lambda_3=\nu$. Denote $u_2=u_{2_1}$, $u_3=u_{3_1}$, and let $\nabla_{u_B}u_C=\Gamma^D_{BC}u_D$, $B, C=1_1, 1_2, 2, 3$, then \eqref{4.1---compatibility-eq1--4} and \eqref{4.1---compatibility-eq2--4} hold.
 If $\lambda_2=\nu$, considering that $\mu$ and $\nu$ are constants, 
we conclude from the Codazzi equation \eqref{Codazzi-td-eq} that
\begin{equation*}
\Gamma_{1_13}^{1_1}=\Gamma_{31_2}^{3}=\Gamma_{31_1}^{3}=\Gamma_{1_23}^{1_1}=
\Gamma_{1_12}^{1_1}=\Gamma_{21_2}^{3}=\Gamma_{31_2}^{2}=0.
\end{equation*}
Applying Gauss equation for $\langle R(u_{1_1}, u_3)u_{1_2}, u_3\rangle$, combining the above equations, we have
$
c+\varepsilon\mu\nu=0.
$

If $\lambda_2=\mu$, then we deduce from the Codazzi equation that 
\begin{equation*}
\Gamma_{1_13}^{1_1}=\Gamma_{1_23}^{1_2}=\Gamma_{1_23}^{1_1}=\Gamma_{31_2}^{3}=\Gamma_{31_1}^{3}=\Gamma_{32}^{3}=\Gamma_{23}^{2}=\Gamma_{23}^{1_1}=\Gamma_{1_23}^{2}=0,
\end{equation*}
and
\begin{equation*}
\Gamma_{21_1}^{3}=\Gamma_{1_12}^{3}.
\end{equation*}
Using Gauss equation for $\langle R(u_{1_1}, u_3)u_{1_2}, u_3\rangle$ and $\langle R(u_{2}, u_3)u_{2}, u_3\rangle$, combining the above two equations, we have
$$
\begin{cases}
\Gamma_{31_2}^{2}\Gamma_{1_12}^{3}=-c\varepsilon_1-\varepsilon\varepsilon_1\mu\nu,\\
-2\Gamma_{31_2}^{2}\Gamma_{1_12}^{3}=-c\varepsilon_1-\varepsilon\varepsilon_1\mu\nu,
\end{cases}
$$
which implies
$
\Gamma_{31_2}^{2}\Gamma_{1_12}^{3}=0,
$
and
$
c+\varepsilon\mu\nu=0.
$

\vskip.2cm
\emph{Subcase (3):\ the shape operator has the form (I\!I\!I).}

\vskip.1cm
Suppose $\lambda_1=\mu$ and $\lambda_2=\nu$.
Let $\nabla_{u_B}u_C=\Gamma^D_{BC}u_D$, $B, C=1_1, 1_2, 2_1, 2_2$.
It follows from Codazzi equation that
$$
\Gamma_{2_11_2}^{2_1}=\Gamma_{1_12_2}^{1_1}=\Gamma_{1_22_2}^{1_1}=0.
$$
Using Gauss equation for $\langle R(u_{1_1}, u_{2_1})u_{1_2}, u_{2_2}\rangle$, combining the above equation, we have
$
c+\varepsilon\mu\nu=0.
$

\vskip.2cm
\emph{Subcase (4):\ the shape operator has the form (I\!V).}

\vskip.1cm
We can suppose $\lambda_1=\mu$ and $\lambda_2=\nu$.
Let $u_2=u_{2_1}$ and $\nabla_{u_B}u_C=\Gamma^D_{BC}u_D$, $B, C=1_1, 1_2, 1_3, 2$,
then \eqref{3.2---compatibility-eq1--4} and \eqref{3.2---compatibility-eq2--4} hold,
and we have from Codazzi equation that
\begin{equation}\label{houbu-eq1-proposition 5.5}
\Gamma_{1_12}^{1_1}=\Gamma_{1_32}^{1_3},\ \ \Gamma_{1_12}^{1_1}+\Gamma_{1_22}^{1_2}+\Gamma_{1_32}^{1_3}=0,
\end{equation}
\begin{equation}\label{houbu-eq2-proposition 5.5}
\Gamma_{21_1}^{2}=\Gamma_{21_2}^{2}=\Gamma_{21_3}^{2}=\Gamma_{1_2 2}^{1_1}=\Gamma_{1_32}^{1_1}=\Gamma_{1_32}^{1_2}=0,
\end{equation}
\begin{equation}\label{houbu-eq3-proposition 5.5}
(\nu-\mu)\Gamma_{1_32}^{1_3}=-\Gamma_{21_3}^{1_2},
\end{equation}
and
\begin{equation}\label{houbu-eq4-proposition 5.5}
(\mu-\nu)\Gamma_{1_11_2}^{2}+\Gamma_{1_11_3}^{2}=(\mu-\nu)\Gamma_{1_21_1}^{2}+\Gamma_{1_21_2}^{2}.
\end{equation}

From Gauss equation for $\langle R(u_{1_1}, u_2)u_{1_3}, u_2\rangle$ and $\langle R(u_{1_3}, u_2)u_{1_1}, u_2\rangle$, using \eqref{houbu-eq2-proposition 5.5}, we have
\begin{equation}\label{houbu-eq5-proposition 5.5}
\begin{cases}
-u_2(\Gamma_{1_13}^{2})+\Gamma_{21_3}^{1_2}\Gamma_{1_11_2}^{2}-\Gamma_{1_12}^{1_1}\Gamma_{1_11_3}^{2}=-c\varepsilon_1-\varepsilon\varepsilon_1\mu\nu,\\
-u_2(\Gamma_{1_31_1}^{2})+\Gamma_{21_3}^{1_2}\Gamma_{1_21_1}^{2}-\Gamma_{1_32}^{1_3}\Gamma_{1_31_1}^{2}=-c\varepsilon_1-\varepsilon\varepsilon_1\mu\nu,
\end{cases}
\end{equation}
which together with \eqref{houbu-eq1-proposition 5.5} implies that
$$
\Gamma_{21_3}^{1_2}=0\ \  \text{or}\ \ \Gamma_{1_11_2}^{2}=\Gamma_{1_21_1}^{2}.
$$
If $\Gamma_{21_3}^{1_2}=0$, then \eqref{houbu-eq3-proposition 5.5} tells us that $\Gamma_{1_32}^{1_3}=0$.
If $\Gamma_{1_11_2}^{2}=\Gamma_{1_21_1}^{2}$, then from \eqref{houbu-eq1-proposition 5.5}, \eqref{houbu-eq3-proposition 5.5} and \eqref{houbu-eq4-proposition 5.5}, we also have
\begin{equation*}
\Gamma_{1_32}^{1_3}=\Gamma_{21_3}^{1_2}=0.
\end{equation*}
So, we derive from \eqref{houbu-eq5-proposition 5.5} that
$
c+\varepsilon\mu\nu=0.
$

\vskip.2cm
\emph{Case 2:\ $\mu$ and $\nu$ are all imaginary.}

\vskip.1cm
We can suppose $\mu=\gamma+\tau\sqrt{-1}$ and $\nu=\gamma-\tau\sqrt{-1}$.
Since $H$ is a nonzero constant, the equations \eqref{proper-eq2} and
\eqref{trA--td---sec2.3} gives that
\begin{equation*}
\begin{cases}
\gamma=\varepsilon H,\\
4(\gamma^2-\tau^2)=\varepsilon \lambda.
\end{cases}
\end{equation*}
So, $\gamma$ and $\tau$ are all constants.
In the following, we treat the possible forms (V\!I\!I) and (V\!I\!I\!I) of $A$, and derive the result that $c+\varepsilon\mu\nu=0$.

\vskip.2cm
\emph{Subcase (1):\ the shape operator has the form (V\!I\!I).}

\vskip.1cm

It follows that $\nu_1=\nu$ and $\tau_1=\tau$.
Let
$$
\begin{aligned}
\nabla_{\bar{u}_{1_a}}\bar{u}_{1_b}=\Gamma_{\bar{1}_a\bar{1}_b}^{\bar{1}_d}\bar{u}_{1_d}+\Gamma_{\bar{1}_a\bar{1}_b}^{\tilde{1}_d}\bar{v}_{1_d},\ \nabla_{\bar{v}_{1_a}}\bar{v}_{1_b}=\Gamma_{\tilde{1}_a\tilde{1}_b}^{\bar{1}_d}\bar{u}_{1_d}+\Gamma_{\tilde{1}_a\tilde{1}_b}^{\tilde{1}_d}\bar{v}_{1_d},\\ \nabla_{\bar{u}_{1_a}}\bar{v}_{1_b}=\Gamma_{\bar{1}_a\tilde{1}_b}^{\bar{1}_d}\bar{u}_{1_d}+\Gamma_{\bar{1}_a\tilde{1}_b}^{\tilde{1}_d}\bar{v}_{1_d},\ \nabla_{\bar{v}_{1_a}}\bar{u}_{1_b}=\Gamma_{\tilde{1}_a\bar{1}_b}^{\bar{1}_d}\bar{u}_{1_d}+\Gamma_{\tilde{1}_a\bar{1}_b}^{\tilde{1}_d}\bar{v}_{1_d},
\end{aligned}
$$
with $a, b=1, 2$,
we get from compatibility condition that
$$
\Gamma_{B\bar{1}_a}^{\bar{1}_b}=-\Gamma_{B\bar{1}_{3-b}}^{\bar{1}_{3-a}},\ \Gamma_{B\tilde{1}_a}^{\tilde{1}_b}=-\Gamma_{B\tilde{1}_{3-b}}^{\tilde{1}_{3-a}},\
\Gamma_{B\bar{1}_a}^{\tilde{1}_b}=\Gamma_{B\tilde{1}_{3-b}}^{\bar{1}_{3-a}}.
$$
Since $\gamma$ and $\tau$ are constants, we deduce from Codazzi equation that
\begin{equation}\label{houbu-eq1-proposition 5.9}
\Gamma_{\bar{1}_2\tilde{1}_2}^{\bar{1}_1}=\Gamma_{\bar{1}_2\bar{1}_2}^{\tilde{1}_1}=\Gamma_{\tilde{1}_2\tilde{1}_2}^{\bar{1}_1}=\Gamma_{\tilde{1}_2\bar{1}_2}^{\tilde{1}_1}=0,
\end{equation}
and
\begin{equation*}
\begin{cases}
-\Gamma_{\tilde{1}_2\tilde{1}_1}^{\tilde{1}_1}+\Gamma_{\tilde{1}_2\bar{1}_1}^{\bar{1}_1}=0,\\
\Gamma_{\bar{1}_2\bar{1}_2}^{\tilde{1}_2}+\Gamma_{\bar{1}_2\tilde{1}_2}^{\bar{1}_2}=-\Gamma_{\tilde{1}_2\tilde{1}_2}^{\tilde{1}_2}+\Gamma_{\tilde{1}_2\bar{1}_2}^{\bar{1}_2},\\
\Gamma_{\bar{1}_1\tilde{1}_2}^{\bar{1}_1}+\Gamma_{\bar{1}_1\bar{1}_2}^{\tilde{1}_1}=\Gamma_{\bar{1}_2\tilde{1}_1}^{\bar{1}_1}+\Gamma_{\bar{1}_2\bar{1}_1}^{\tilde{1}_1},\\
-\Gamma_{\bar{1}_2\tilde{1}_1}^{\tilde{1}_1}+\Gamma_{\bar{1}_2\bar{1}_1}^{\bar{1}_1}=0,\\
\Gamma_{\bar{1}_2\bar{1}_2}^{\bar{1}_2}-\Gamma_{\bar{1}_2\tilde{1}_2}^{\tilde{1}_2}=-\Gamma_{\tilde{1}_2\tilde{1}_2}^{\bar{1}_2}-\Gamma_{\tilde{1}_2\bar{1}_2}^{\tilde{1}_2},\\
\Gamma_{\tilde{1}_1\bar{1}_2}^{\tilde{1}_1}+\Gamma_{\tilde{1}_1\tilde{1}_2}^{\bar{1}_1}=\Gamma_{\tilde{1}_2\bar{1}_1}^{\tilde{1}_1}+\Gamma_{\tilde{1}_2\tilde{1}_1}^{\bar{1}_1},
\end{cases}
\end{equation*}
which reduces to
\begin{equation}\label{houbu-eq2-proposition 5.9}
\Gamma_{\bar{1}_1\bar{1}_2}^{\tilde{1}_1}=\Gamma_{\bar{1}_1\tilde{1}_2}^{\bar{1}_1}=\Gamma_{\tilde{1}_1\bar{1}_2}^{\tilde{1}_1}=\Gamma_{\tilde{1}_1\tilde{1}_2}^{\bar{1}_1}=0.
\end{equation}
Using Gauss equation for $\langle R(\bar{u}_{1_1}, \bar{v}_{1_1})\bar{u}_{1_2}, \bar{v}_{1_2}\rangle$, combining \eqref{houbu-eq1-proposition 5.9} and \eqref{houbu-eq2-proposition 5.9}, we have
$$
c+\varepsilon(\gamma^2+\tau^2)=0,
$$
i.e. $c+\varepsilon\mu\nu=0$.

\vskip.2cm
\emph{Subcase (2):\ the shape operator has the form (V\!I\!I\!I).}

\vskip.1cm
In this case, we have $\gamma_1=\gamma_2=\gamma$, $\tau_1=\tau_2=\tau$.
Denote $\bar{u}_1=\bar{u}_{1_1}$, $\bar{u}_2=\bar{u}_{2_1}$ and $\bar{v}_1=\bar{v}_{1_1}$ and $\bar{v}_2=\bar{v}_{2_1}$.
Let
$$
\begin{aligned}
\nabla_{\bar{u}_{a}}\bar{u}_{b}=\Gamma_{\bar{a}\bar{b}}^{\bar{d}}\bar{u}_{d}+\Gamma_{\bar{a}\bar{b}}^{\tilde{d}}\bar{v}_{d},\ \nabla_{\bar{v}_{a}}\bar{v}_{b}=\Gamma_{\tilde{a}\tilde{b}}^{\bar{d}}\bar{u}_{d}+\Gamma_{\tilde{a}\tilde{b}}^{\tilde{d}}\bar{v}_{d},\\ \nabla_{\bar{u}_{a}}\bar{v}_{b}=\Gamma_{\bar{a}\tilde{b}}^{\bar{d}}\bar{u}_{d}+\Gamma_{\bar{a}\tilde{b}}^{\tilde{d}}\bar{v}_{d},\ \nabla_{\bar{v}_{a}}\bar{u}_{b}=\Gamma_{\tilde{a}\bar{b}}^{\bar{d}}\bar{u}_{d}+\Gamma_{\tilde{a}\bar{b}}^{\tilde{d}}\bar{v}_{d},
\end{aligned}
$$
with $a, b=1, 2$,
we obtain from compatibility condition that
$$
\Gamma_{B\bar{a}}^{\bar{b}}=-\Gamma_{B\bar{b}}^{\bar{a}},\ \Gamma_{B\tilde{a}}^{\tilde{b}}=-\Gamma_{B\tilde{b}}^{\tilde{a}},\
\Gamma_{B\bar{a}}^{\tilde{b}}=\Gamma_{B\tilde{b}}^{\bar{a}},\ \ a, b=1, 2.
$$
The Codazzi equation gives that
\begin{equation*}
\Gamma_{\bar{1}\tilde{1}}^{\bar{1}}=\Gamma_{\tilde{1}\bar{1}}^{\tilde{1}}=\Gamma_{\bar{2}\bar{1}}^{\tilde{1}}=\Gamma_{\tilde{2}\bar{1}}^{\tilde{1}}=0,
\end{equation*}
and
\begin{equation*}
\Gamma_{\bar{1}\tilde{2}}^{\bar{1}}+\Gamma_{\bar{1}\bar{2}}^{\tilde{1}}=\Gamma_{\bar{1}\bar{2}}^{\bar{1}}-\Gamma_{\bar{1}\tilde{2}}^{\tilde{1}}=\Gamma_{\tilde{1}\tilde{2}}^{\bar{1}}+\Gamma_{\tilde{1}\bar{2}}^{\tilde{1}}=\Gamma_{\tilde{1}{\tilde{2}}}^{\tilde{1}}-\Gamma_{\tilde{1}\bar{2}}^{\bar{1}}=0.
\end{equation*}
Applying Gauss equation for $\langle R(\bar{u}_{1}, \bar{v}_{1})\bar{u}_{1}, \bar{v}_{1}\rangle$, combining the above two equations, we conclude
$
c+\varepsilon(\gamma^2+\tau^2)=0,
$
i.e. $c+\varepsilon\mu\nu=0$.
$\hfill\square$

\vskip.2cm
\noindent
{\bf Remark 4.4}\quad The result of Lemma 4.3 that $c+\varepsilon\mu\nu=0$ coincides with the basic identity of Cartan in [14, Theorem 2.9] for the isoparametric hypersurface $M^n_r$ of $N^{n+1}_s(c)(c=0,-1,1)$ with two distinct principal curvatures $\mu$ and $\nu$, and algebraic and geometric multiplicities of $\mu$ or $\nu$ coincide.

\vskip.2cm

By use of Lemma 4.3, we can give the value of the mean curvature $H$, according to the principal curvatures are real or imaginary. The following Lemma 4.5 will be used for the case that the principal curvatures are real.

\vskip.2cm
\noindent
{\bf Lemma 4.5}\quad \emph{Let $c$, $\lambda$, $l$, $\mu$ and $\nu$ are real constants, and the equations
$$
\begin{cases}
\mu\nu=-c\varepsilon,\\
l\mu^2+(4-l)\nu^2=\varepsilon\lambda
\end{cases}
$$
hold, where $\varepsilon=\pm 1$. Then, $\varepsilon\lambda=-4c\varepsilon\geq 0$ is the necessary condition such that $\mu=\nu$.}

\vskip.1cm
{\bf Proof}\quad If $\mu=\nu$, then
$$
-c\varepsilon=\mu^2\geq 0,
$$ and
$$
\varepsilon\lambda=-4c\varepsilon.
$$
$\hfill\square$

\vskip.2cm
\noindent
{\bf Theorem 4.6}\quad \emph{Let $M^4_r$ be a nondegenerate hypersurface of $N^{5}_s(c)$ with proper mean curvature vector field, satisfying $\Delta\vec{H}=\lambda\vec{H}$, $\vec{\xi}$ be a unit normal vector field to $M^4_r$, with $\varepsilon=\langle\vec\xi,\vec\xi\rangle=\pm1$. Suppose that $M^4_r$ has two distinct real principal curvatures $\mu$ and $\nu$, and $l=1,\ 2\ \text{or}\ 3$ is the multiply of $\mu$, then
\begin{itemize}
  \item When $c=0$, we have
  \begin{itemize}
    \item if $\varepsilon \lambda\leq 0$, then $H=0$;
    \item if $\varepsilon \lambda>0$, then $H=0$, or
$
H^2=\frac{l\varepsilon\lambda}{16}
$, or $
H^2=\frac{(4-l)\varepsilon\lambda}{16}
$;
  \end{itemize}
  \item When $c\neq 0$, we have
  \begin{itemize}
    \item if $\varepsilon \lambda<2\sqrt{l(4-l)}|c|$, then $H=0$;
    \item if $\varepsilon \lambda\geq2\sqrt{l(4-l)}|c|$, then $H=0$ or
$$
H^2=\frac{1}{16}[2\varepsilon\lambda\pm(l-2)\sqrt{\lambda^2-4l(4-l)c^2}-2l(4-l)c\varepsilon].
$$
Specially, if $\varepsilon\lambda=-4c\varepsilon>0$, then $H=0$ or $H^2=-\frac{3}{4}c\varepsilon$.
  \end{itemize}
\end{itemize}
}

{\bf Proof}
\quad (i)\ When $c=0$, suppose $H\neq 0$, then \eqref{trA-trA2-td-tworeal-bu1-sec2.3} holds, and
$
\varepsilon \lambda>0.
$
It follows from Lemma 4.3 that $\mu\nu=0$, which together with the second equation in \eqref{trA-trA2-td-tworeal-bu1-sec2.3} gives that
$$
\mu=0,\ \nu^2=\frac{\varepsilon\lambda}{4-l},\ \ \text{or}\ \mu^2=\frac{\varepsilon\lambda}{l},\ \nu=0.
$$
Since the above equation, we obtain from the first equation in \eqref{trA-trA2-td-tworeal-bu1-sec2.3} that
$$
H^2=\frac{l\varepsilon\lambda}{16},\ \ \text{or}\ H^2=\frac{(4-l)\varepsilon\lambda}{16}.
$$

(ii)\ For the case $c\neq 0$, suppose $H\neq 0$, then \eqref{trA-trA2-td-tworeal-bu1-sec2.3} holds, and $
\mu\nu=-c\varepsilon$, obtained from Lemma 4.3.
And then, we have
\begin{equation*}
\begin{cases}
\varepsilon \lambda=\sqrt{l}\mu^2+\sqrt{4-l}\nu^2>0,\\
\varepsilon \lambda+2\sqrt{l(4-l)}c\varepsilon=(\sqrt{l}\mu-\sqrt{4-l}\nu)^2\geq 0,\\
\varepsilon \lambda-2\sqrt{l(4-l)}c\varepsilon=(\sqrt{l}\mu+\sqrt{4-l}\nu)^2\geq 0,
\end{cases}
\end{equation*}
which implies $\varepsilon\lambda\geq2\sqrt{l(4-l)}|c|$. So, if $\varepsilon\lambda<2\sqrt{l(4-l)}|c|$, then $H=0$.

When $\varepsilon \lambda\geq2\sqrt{l(4-l)}|c|$, if $H\neq 0$, then we calculate $\mu^2$ and $\nu^2$ from $\mu\nu=-4c\varepsilon$ and the second equation of
\eqref{trA-trA2-td-tworeal-bu1-sec2.3} that
\begin{equation*}
\mu^2=\frac{\varepsilon\lambda\pm\sqrt{\lambda^2-4l(4-l)c^2}}{2l},\ \ \nu^2=\frac{\varepsilon\lambda\mp\sqrt{\lambda^2-4l(4-l)c^2}}{2(4-l)}.
\end{equation*}
Together with the first equation of \eqref{trA-trA2-td-tworeal-bu1-sec2.3} and the above equation, we have
\begin{equation*}
H^2=\frac{2\varepsilon\lambda\pm(l-2)\sqrt{\lambda^2-4l(4-l)c^2}-2l(4-l)c\varepsilon}{16}.
\end{equation*}

However, some values of $\mu^2$, $\nu^2$ and $H^2$ in the above maybe contradict with the condition that $\mu\neq\nu$. From Lemma 4.5, considering $\varepsilon \lambda\neq 0$, we know $\varepsilon\lambda=-4c\varepsilon>0$ is the necessary condition of $\mu=\nu$ (a contradiction).
It's necessary for us to discuss the values of $\mu^2$, $\nu^2$ and $H^2$ for the case $\varepsilon\lambda=-4c\varepsilon>0$, which is in the range $\varepsilon \lambda\geq2\sqrt{l(4-l)}|c|$.
By calculation, when $\varepsilon\lambda=-4c\varepsilon>0$, we find

\vskip.1cm
(1)\ $\mu^2=\nu^2=H^2=-c\varepsilon$ for $l=2$;

\vskip.1cm
(2)\ $\mu^2=-3c\varepsilon$, $\nu^2=-\frac{c\varepsilon}{3}$, $H^2=-\frac{3}{4}c\varepsilon$, or $\mu^2=\nu^2=H^2=-c\varepsilon$, for $l=1$;

\vskip.1cm
(3)\ $\nu^2=-3c\varepsilon$, $\mu^2=-\frac{c\varepsilon}{3}$, $H^2=-\frac{3}{4}c\varepsilon$, or $\mu^2=\nu^2=H^2=-c\varepsilon$, for $l=3$.\\
As $\mu\nu>0$, $\mu^2=\nu^2$ is equivalent to $\mu=\nu$ (a contradiction). Thus, we conclude that $H^2=-\frac{3}{4}c\varepsilon$ for the case $\varepsilon\lambda=-4c\varepsilon>0$.

In conclusion, when $c\neq 0$ and $\varepsilon \lambda\geq2\sqrt{l(4-l)}|c|$, we have $H=0$ or
\begin{equation}\label{H2-sec4.2}
H^2=\frac{2\varepsilon\lambda\pm(l-2)\sqrt{\lambda^2-4l(4-l)c^2}-2l(4-l)c\varepsilon}{16}.
\end{equation}
Specially, if $\varepsilon\lambda=-4c\varepsilon>0$, then $H=0$ or $H^2=-\frac{3}{4}c\varepsilon$.

By the way, we emphasize that the right side values of \eqref{H2-sec4.2} are all greater than or equal to zero.
When $l=2$, it is obviously.
When $l=1, 3$, since $\varepsilon \lambda\geq2\sqrt{3}|c|$, it follows that
$$
\varepsilon\lambda-3|c|\geq(2\sqrt{3}-3)|c|\geq 0,
$$
which together with
$$
\frac{(2\varepsilon\lambda-6c\varepsilon)^2}{\lambda^2-12c^2}=1+\frac{3(\lambda-4c)^2}{\lambda^2-12c^2}\geq 1
$$
gives that $2\varepsilon\lambda-6c\varepsilon\geq\sqrt{\lambda^2-12c^2}$.
So, when $l=1, 3$, the right side values
$$
\frac{2\varepsilon\lambda\pm\sqrt{\lambda^2-12c^2}-6c\varepsilon}{16}\geq 0.
$$

$\hfill\square$

\vskip.2cm
\noindent
{\bf Remark 4.7}\quad When $c=0$ and the algebraic and geometric multiplicities of $\mu$ or $\nu$ coincide, we has gotten in
\cite{Liu 2017} the values of $H^2$, which is agree with the result of Theorem 4.6 for $c=0$.

\vskip.2cm
\noindent
{\bf Theorem 4.8}\quad \emph{Let $M^4_r$ be a nondegenerate hypersurface of $N^{5}_s(c)$ with proper mean curvature vector field, satisfying $\Delta\vec{H}=\lambda\vec{H}$, $\vec{\xi}$ be a unit normal vector field to $M^4_r$, with $\varepsilon=\langle\vec\xi,\vec\xi\rangle=\pm1$. Suppose that $M^4_r$ has two imaginary principal curvatures, we have
\begin{itemize}
  \item If $c\varepsilon\geq0$, then $H= 0$;
  \item If $c\varepsilon<0$ and $|\varepsilon\lambda|\geq -4c\varepsilon$, then $H=0$;
  \item If $c\varepsilon<0$ and $|\varepsilon\lambda|<-4c\varepsilon$, then $H=0$ or
$
H^2=\frac{\varepsilon\lambda-4c\varepsilon}{8}
$.
\end{itemize}
}

{\bf Proof}\quad Suppose that $M^4_r$ has two imaginary principal curvatures $\gamma+\tau i$ and $\gamma-\tau i$, we get from Lemma 4.3 and its proof that if $H\neq 0$, then $\gamma=\varepsilon H$ and
\begin{equation}\label{houbu-eq1-theorem 5.16}
\begin{cases}
4(\gamma^2-\tau^2)=\varepsilon \lambda,\\
\gamma^2+\tau^2=-c\varepsilon.
\end{cases}
\end{equation}

When $c\varepsilon\geq 0$, we assume that $H\neq 0$, then $\gamma^2+\tau^2=-c\varepsilon\leq 0$, a contradiction.

When $c\varepsilon<0$ and $|\varepsilon\lambda|\geq -4c\varepsilon$, we suppose $H\neq 0$, then the equations
\eqref{houbu-eq1-theorem 5.16} implies that
$$
\begin{cases}
-4c\varepsilon-\varepsilon \lambda=8\tau^2>0,\\
\varepsilon \lambda-4c\varepsilon=8\nu^2> 0,
\end{cases}
$$
which contradicts with $|\varepsilon\lambda|\geq -4c\varepsilon$.

For the case $c\varepsilon<0$ and $|\varepsilon\lambda|<-4c\varepsilon$, if $H\neq 0$, we also have $8\gamma^2=\varepsilon \lambda-4c\varepsilon$, which combining $\gamma=\varepsilon H$ tells us that
$$
H^2=\frac{\varepsilon\lambda-4c\varepsilon}{8}.
$$
$\hfill\square$

\vskip.2cm
\noindent
{\bf Corollary 4.9}\quad \emph{Let $M^4_r$ be a nondegenerate hypersurface of $\mathbb{E}^{5}_s$ with proper mean curvature vector field. Suppose that $M^4_r$ has two principal curvatures, which are imaginary,
then $M^4_r$ is minimal.}

\vskip.2cm
\noindent
{\bf Remark 4.10}\quad In \cite{Liu 2017}, we also estimated the value of the mean curvature $H$ for the hypersurface $M^4_r$ in $\mathbb{E}^5_s$ satisfying $\Delta\vec{H}=\lambda\vec{H}$ and with two distinct imaginary principal curvatures, but just gave a value range that $H=0$ or $H^2>\frac{\varepsilon\lambda}{4}$. Clearly,
the result $H=0$ of Corollary 4.9 improves the result in \cite{Liu 2017} to the greatest extent.

\vskip.2cm

In non-flat space form $N^{5}_s(c)\ (c\neq 0)$, under the assumption that the hypersurface $M^4_r$ has two imaginary principal curvatures, we have from Theorem 4.8 that the mean curvature $H$ is not necessarily zero. However, if the hypersurface $M^4_r$ is biharmonic, i.e. $\lambda=4c$, then Theorem 4.8 yields that the mean curvature is zero, which give a partial affirmative answer to Chen's conjecture.

\vskip.2cm
\noindent
{\bf Corollary 4.11}\quad \emph{Let $M^4_r$ be a biharmonic hypersurface of $N^{5}_s(c)$. Suppose that $M^4_r$ has two principal curvatures, which are imaginary,
then $M^4_r$ is minimal.}

\subsection{When $M^4_r$ has not imaginary principal curvatures}

\vskip.2cm
\noindent
{\bf Theorem 4.12}\quad \emph{Let $M^4_r$ be a nondegenerate hypersurface of $N^{5}_s(c)$ satisfying $\Delta\vec{H}=\lambda\vec{H}$ ($\lambda$ a constant), and $\vec{\xi}$ be a unit normal vector field to $M^4_r$, with $\varepsilon=\langle\vec\xi,\vec\xi\rangle=\pm1$. Suppose that $M^4_r$ has not imaginary principal curvatures, we have
\begin{itemize}
  \item If $\varepsilon \lambda\leq 0$, then $H=0$;
  \item If $\varepsilon \lambda> 0$, then $H^2\leq\frac{\varepsilon \lambda}{4}$, the equal sign holds if and only if the principal curvatures are all equal and nonzero.
\end{itemize}
}

{\bf Proof}\quad Since $M^4_r$ has not imaginary principal curvatures, the equations \eqref{trA--td---sec2.3} reduce to
\begin{equation*}
\text{tr}A^2=\sum_{i=1}^m \alpha_i\lambda_i^2,\ \
4\varepsilon H=\sum_{i=1}^m \alpha_i\lambda_i.
\end{equation*}
Note that $\sum_{i=1}^m \alpha_i=4$.
Because of that the mean curvature $H$ is a constant, combining the above equations, we have from \eqref{proper-eq2} that
$H=0$ or $\sum_{i=1}^m \alpha_i\lambda_i^2=\varepsilon \lambda$.

When $\varepsilon\lambda\leq 0$, we assume that $H\neq 0$, then
$$
\sum_{i=1}^m \alpha_i\lambda_i^2=\varepsilon \lambda\leq 0,
$$
which implies $\lambda_i=0$, $i=1, \cdots, m$. And then, $4\varepsilon H=\sum_{i=1}^m \alpha_i\lambda_i=0$, a contradiction.

For the case $\varepsilon\lambda>0$, if $H\neq 0$,
using Cauchy inequality, we know
\begin{equation}\label{bu-zh1--sec2.3}
(\sum_{i=1}^m \alpha_i\lambda_i)^2\leq 4\sum_{i=1}^m \alpha_i\lambda_i^2,
\end{equation}
where the equal sign holds if and only if $\lambda_1=\lambda_2=\cdots=\lambda_m\neq 0$.
The equation \eqref{bu-zh1--sec2.3} yields that
$$
0< H^2\leq\frac{\varepsilon \lambda}{4}.
$$
The equal sign of the above equation is true if and only if the principal curvatures are all equal and nonzero.
$\hfill\square$

\vskip.2cm
\noindent
{\bf Corollary 4.13}\quad \emph{Let $M^4_r$ be a biharmonic hypersurface of $N^{5}_s(c)$, $\vec{\xi}$ be a unit normal vector field to $M^4_r$, with $\varepsilon=\langle\vec\xi,\vec\xi\rangle=\pm1$. If $\varepsilon c\leq 0$ and $M^4_r$ has not imaginary principal curvatures,
then $M^4_r$ is minimal.}

\vskip.2cm
\noindent
{\bf Remark 4.14}\quad
We claim that the value
\small{$
\frac{1}{16}[2\varepsilon\lambda\pm(l-2)\sqrt{\lambda^2-4l(4-l)c^2}-2l(4-l)c\varepsilon]
$}
of $H^2$ for the case $c\neq 0$, $\varepsilon \lambda\geq2\sqrt{l(4-l)}|c|$ and $\varepsilon\lambda\neq -4c\varepsilon$ in Theorem 4.6 is in the range $H^2<\frac{\varepsilon \lambda}{4}$.
In fact, in this case, we can easily check that
$
\varepsilon \lambda+l(4-l)c\varepsilon\geq 0$
and
\begin{equation*}
 (2\varepsilon \lambda+2l(4-l)c\varepsilon)^2-(l-2)^2(\lambda^2-4l(4-l)c^2)=l(4-l)(\lambda+4c)^2>0.
\end{equation*}
Thus, it follows that
 $$
 2\varepsilon \lambda+2l(4-l)c\varepsilon>\pm(l-2)\sqrt{\lambda^2-4l(4-l)c^2},
  $$
which deduce our claim.

\vskip.2cm
\noindent
{\bf Remark 4.15}\quad When $M^4_r$ has imaginary principal curvature, we have from \eqref{trA--td---sec2.3} that
$$
\text{tr}A^2=\sum_{i=1}^t\alpha_i\lambda_i^2+2\sum_{j=t+1}^m\beta_j(\nu_j^2-\tau_j^2),\ \ t< m.
$$
Since there exists negative terms in the right side of the above equation, we can not get the range value of $H$ by using Cauchy inequality, for the case that the number of distinct
principal curvatures is larger than 2.

\vskip.3cm\noindent
{\bf Acknowledgements}\quad This work was supported by the National Natural Science Foundation
of China (Nos. 12161078, 11761061), the Foundation for Distinguished Young Scholars of Gansu Province (No. 20JR5RA515),
and the Project of Northwest Normal University (No. NWNU-LKQN2019-23).

\bibliographystyle{model1a-num-names}

\vskip.3cm\noindent
{\bf References}
\begin{thebibliography}{99}
\bibitem{Arvanitoyeorgos 2007}
A. Arvanitoyeorgos, F. Defever, G. Kaimakamis,
Hypersurfaces of $E^4_s$ with proper mean curvature vector,
J. Math. Soc. Japan
{\bf 59}:3 (2007), 797--809.

\bibitem{Arvanitoyeorgos 2009}
A. Arvanitoyeorgos, G. Kaimakamis, M. Magid,
Lorentz hypersurfaces in $E^4_1$ satisfying $\Delta\vec{H}=\alpha\vec{H}$,
Illinois J. Math.
{\bf 53}:2 (2009), 581--590.

\bibitem{Arvanitoyeorgos 2013}
A. Arvanitoyeorgos, G. Kaimakamis,
Hypersurfaces of type $M^3_2$ in $E^4_2$ with proper mean curvature vector,
J. Geom. Phys.
{\bf 63} (2013), 99--106.

\bibitem{Chen 1988-1}
B.-Y. Chen,
Null 2-type surfaces in $\mathbb{E}^3$ are circular cylinders,
Kodai Math. J.
{\bf 11}:2 (1988), 295--299.

\bibitem{Chen 1991}
B.-Y. Chen,
Some open problems and conjectures on submanifolds of finite type,
Soochow J. Math.
{\bf 17}:2 (1991), 169--188.

\bibitem{Defever 1997}
F. Defever,
Hypersurfaces of $\mathbb{E}^4$ satisfying $\Delta\vec{H}=\lambda\vec{H}$,
Michigan Math. J.
{\bf 44}:2 (1997), 355--363.

\bibitem{Du 2017}
L. Du, J. Zhang, X. Xie,
Hypersurfaces satisfying $\tau_2(\phi)=\eta\tau(\phi)$ in pseudo-Riemannian space forms,
Math. Phys. Anal. Geom.
{\bf 20}:2 (2017), 17. 

\bibitem{Du 2022}
L. Du, J. Ren,
On $\eta$-biharmonic hypersurfaces in pseudo-Riemannian space forms,
Math. Slovaca,
{\bf 72} (2022), 1259--1272.

\bibitem{Du 2023}
L. Du,
On $\eta$-biharmonic hypersurfaces with constant scalar curvature in higher dimensional pseudo-Riemannian space forms,
J. Math. Anal. Appl.
{\bf 518} (2023), 122670.

\bibitem{Ferrandez 1992-1}
A. Ferr\'{a}ndez, P. Lucas,
On surfaces in the 3-dimensional Lorentz-Minkowski space,
Pacific J. Math.
{\bf 152}:1 (1992), 93--100.


\bibitem{Fu 2021}
Y. Fu, M. C. Hong, X. Zhan,
On Chen's biharmonic conjecture for hypersurfaces in $\mathbb{R}^5$,
Adv. Math.
{\bf 383} (2021), 107697. 

\bibitem{Fu 2021-2}
Y. Fu, X. Zhan,
Hypersurfaces satisfying $\Delta\vec{H}=\lambda\vec{H}$ in $\mathbb{E}^5$,
J. Math. Anal. Appl.
{\bf 503}:2 (2021), 125337. 

\bibitem{Inoguchi 2007}
J.-I. Inoguchi,
Biminimal submanifolds in contact 3-manifolds,
Balkan J. Geom. Appl.
{\bf 12}:1 (2007), 56--67.

\bibitem{Liu 2014}
J.-C. Liu, C. Yang,
Hypersurfaces in $\mathbb{E}^{n+1}_s$ satisfying
$\Delta \vec{H}=\lambda \vec{H}$ with at most three distinct principal curvatures,
J. Math. Anal. Appl.
{\bf 419} (2014), 562--573.


\bibitem{Liu 2016}
J.-C. Liu, C. Yang,
Lorentz hypersurfaces in $\mathbb{E}^{n+1}_1$ satisfying
$\Delta \vec{H}=\lambda \vec{H}$ with at most three distinct principal curvatures,
J. Math. Anal. Appl.
{\bf 434} (2016), 222--240.

\bibitem{Liu 2017}
J.-C. Liu, C. Yang,
Hypersurfaces in $\mathbb{E}^{n+1}_s$ satisfying $\Delta\vec{H}=\lambda\vec{H}$ with at most two
distinct principal curvatures,
J. Math. Anal. Appl.
{\bf 451} (2017), 14--33.

\bibitem{Neill 1983}
B. O'Neill,
Semi-Riemannian Geometry: with Applications to Relativity,
Pure Appl. Math. 103, Academic Press, New York,
1983.


\bibitem{Hasanis 1995}
T. Hasanis, T. Vlachos,
Hypersurfaces with constant scalar curvature and constant mean curvature,
Ann. Global Anal. Geom.
{\bf 13}:1 (1995), 69--77.

\end{thebibliography}

\end{document}